\documentclass{article}
\usepackage[utf8]{inputenc}
\usepackage[spanish,english]{babel}
\usepackage{amsmath,amssymb,amsthm,amsfonts}
\usepackage{geometry}
\usepackage{hyperref}
\usepackage{fancyhdr}
\usepackage{titlesec}
\usepackage{listings}
\usepackage{graphicx,graphics}
\usepackage{multicol}
\usepackage{multirow}
\usepackage{color}
\usepackage{float} 
\usepackage{subfig}
\usepackage[figuresright]{rotating}
\usepackage{enumerate}
\usepackage{anysize} 
\usepackage{url}
\title{Procesos Regenerativos y de Renovación: Revisión\\
\small{Renewal and Regenerative Processes: Review}}
\author{Carlos E. Martínez-Rodríguez\\
	Academia de Matem\'aticas\\
	Universidad Aut\'onoma de la Ciudad de M\'exico\\
	Iztapalapa,  Ermita Iztapalapa 4163\\
	\texttt{carlos.martinez@uacm.edu.mx} }
\date{}
\newtheorem{Def}{Definición}[section]
\newtheorem{Ejem}{Ejemplo}[section]
\newtheorem{Teo}{Teorema}[section]
\newtheorem{Note}{Nota}[section]
\newtheorem{Prop}{Proposición}[section]
\newtheorem{Cor}{Corolario}[section]
\newtheorem{Lema}{Lema}[section]
\newtheorem{Sup}{Supuestos}[section]

\newcommand{\nat}{\mathbb{N}}

\newcommand{\rea}{\mathbb{R}}
\newcommand{\Eb}{\mathbf{E}}
\newcommand{\esp}{\mathbb{E}}
\newcommand{\prob}{\mathbb{P}}
\newcommand{\indora}{\mbox{$1$\hspace{-0.8ex}$1$}}
\newcommand{\ER}{\left(E,\mathcal{E}\right)}
\newcommand{\KM}{\left(P_{s,t}\right)}
\newcommand{\PE}{\left(X_{t}\right)_{t\in I}}

\renewcommand{\abstractname}{Resumen}
\numberwithin{equation}{section}

\begin{document}
\maketitle

\begin{abstract}
En este documento se presenta una recopilación de resultados relacionados con la teoría de procesos estocásticos, con un enfoque específico en procesos de Markov, procesos regenerativos, procesos de renovación y procesos estacionarios. La relevancia de estos temas reside en la capacidad de identificar puntos de regeneración y las condiciones necesarias para la garantizar la estacionareidad del proceso. El estudio inicia con una revisión de cadenas de Markov y prosigue con el análisis de procesos que cumplen con la propiedad fuerte de Markov. Posteriormente, se profundiza en los procesos de renovación, procesos regenerativos y, finalmente, en los procesos regenerativos estacionarios, destacando los resultados presentados por Thorisson \cite{Thorisson}. Este trabajo no tiene la intención de ser exhaustivo, sino de proporcionar una base sólida que permita profundizar en el conocimiento de estos procesos, dado su amplio espectro de aplicaciones en criptografía \cite{Portman}, teoría de colas \cite{Morozov} y métodos de Monte Carlo \cite{Xu}. Asimismo, se subraya la importancia de los procesos de tipo Poisson debido a sus numerosas aplicaciones (ver \cite{Douc}).
\end{abstract}

\begin{otherlanguage}{english}
\renewcommand{\abstractname}{Abstract} 
\begin{abstract}
This document presents a compilation of results related to the theory of stochastic processes, with a specific focus on Markov processes, regenerative processes, renewal processes, and stationary processes. The relevance of these topics lies in the ability to identify regeneration points and the necessary conditions to ensure the stationarity of the process. The study begins with a review of Markov chains and continues with the analysis of processes that satisfy the strong Markov property. Subsequently, it delves into renewal processes, regenerative processes, and finally, stationary regenerative processes, highlighting the results presented by Thorisson \cite{Thorisson}. This work is not intended to be exhaustive but aims to provide a solid foundation for further deepening the knowledge of these processes, given their broad range of applications in cryptography \cite{Portman}, queueing theory \cite{Morozov}, and Monte Carlo methods \cite{Xu}. Additionally, the importance of Poisson-type processes is emphasized due to their numerous applications (see \cite{Douc}).
\end{abstract}
\end{otherlanguage}


\section*{Introducción}
El estudio de las cadenas de Markov es fundamental para comprender las condiciones bajo las cuales un proceso estocástico puede regenerarse, así como para determinar la existencia de tiempos de regeneración. La extensión de estos conceptos a teorías de colas y sistemas de visitas cíclicas requiere un conocimiento profundo de la teoría subyacente. Este análisis naturalmente lleva a explorar temas más complejos, tales como los procesos regenerativos y de renovación. En este trabajo se realiza una revisión de los conceptos esenciales para iniciar el estudio de los procesos regenerativos estacionarios. La revisión de estos temas se llevó a cabo en su momento bajo la supervisión del Dr. Raúl Montes de Oca Machorro y la Dra. Patricia Saavedra Barrera, cuyas oportunas y valiosas sugerencias y comentarios fueron fundamentales para desarrollar el estudio de este tipo de procesos estocásticos. Es importante destacar que las aplicaciones de estos resultados en la teoría de colas tienen un impacto significativo en problemas contemporáneos. A pesar de los avances logrados, aún existen preguntas sin resolver sobre su aplicación y generalización en la teoría de colas. Este trabajo no pretende ser un estudio exhaustivo sobre el tema, sino más bien proporcionar los elementos necesarios para introducirse en el estudio de estos procesos. El documento está organizado de la siguiente manera: en la primera sección se realiza una revisión de las cadenas de Markov y de los procesos de Markov. En la segunda sección se aborda un estudio inicial de los procesos de renovación y de los procesos regenerativos, junto con sus propiedades y el teorema principal de renovación. La tercera sección profundiza en los procesos regenerativos incluidos en \cite{Thorisson}, para los cuales es necesario revisar procesos más generales. Finalmente, en la última sección se presentan una serie de consideraciones respecto al contenido de este trabajo.

\begin{otherlanguage}{english}
\renewcommand{\abstractname}{Abstract} 
\section*{Introduction}

The study of Markov chains is fundamental for understanding the conditions under which a stochastic process can regenerate, as well as for determining the existence of regeneration times. Extending these concepts to queueing theories and cyclic visit systems requires a deep understanding of the underlying theory. This analysis naturally leads to exploring more complex topics such as regenerative and renewal processes. This work provides a review of the essential concepts for initiating the study of stationary regenerative processes. The review of these topics was conducted under the supervision of Dr. Raúl Montes de Oca Machorro and Dr. Patricia Saavedra Barrera, whose timely and valuable suggestions and comments were fundamental in developing the study of these types of stochastic processes. It is important to highlight that the applications of these results in queueing theory have a significant impact on contemporary problems. It is important to highlight that the applications of these results in queueing theory have a significant impact on contemporary problems. Despite the advances made, there are still unresolved questions regarding their application and generalization in queueing theory. This work does not aim to be an exhaustive study of the topic but rather to provide the necessary elements to introduce the study of these processes. The document is organized as follows: the first section provides a review of Markov chains and Markov processes. The second section addresses an initial study of renewal processes and regenerative processes, along with their properties and the main renewal theorem. The third section delves into the regenerative processes included in \cite{Thorisson}, for which it is necessary to review more general processes. Finally, the last section presents a series of considerations regarding the content of this work.

\end{otherlanguage}

\section{Procesos Estoc\'asticos}\label{Procesos.Estocasticos}

\begin{Def}\index{Conjunto Medible}
Sea $X$ un conjunto y $\mathcal{F}$ una $\sigma$-\'algebra de subconjuntos de $X$, la pareja $\left(X,\mathcal{F}\right)$ es llamado espacio medible. Un subconjunto $A$ de $X$ es llamado medible, o medible con respecto a $\mathcal{F}$, si $A\in\mathcal{F}$.
\end{Def}

\begin{Def}\index{Medida $\sigma$-finita}
Sea $\left(X,\mathcal{F},\mu\right)$ espacio de medida. Se dice que la medida $\mu$ es $\sigma$-finita si se puede escribir $X=\bigcup_{n\geq1}X_{n}$ con $X_{n}\in\mathcal{F}$ y $\mu\left(X_{n}\right)<\infty$.
\end{Def}

\begin{Def}\label{Cto.Borel}\index{Conjunto de Borel}
Sea $X$ un espacio topológico. El álgebra de Borel en $X$, denotada por $\mathcal{B}(X)$, es la $\sigma$-álgebra generada por la colección de todos los conjuntos abiertos de $X$. Es decir, $\mathcal{B}(X)$ es la colección más pequeña de subconjuntos de $X$ que contiene todos los conjuntos abiertos y es cerrada bajo la unión numerable, la intersección numerable y el complemento.
\end{Def}

\begin{Def}\label{Funcion.Medible}\index{Funci\'on Medible}
Una funci\'on $f:X\rightarrow\rea$, es medible si para cualquier n\'umero real $\alpha$ el conjunto \[\left\{x\in X:f\left(x\right)>\alpha\right\},\] pertenece a $X$. Equivalentemente, se dice que $f$ es medible si \[f^{-1}\left(\left(\alpha,\infty\right)\right)=\left\{x\in X:f\left(x\right)>\alpha\right\}\in\mathcal{F}.\]
\end{Def}

\begin{Def}\label{Def.Cilindros}\index{Cilindro}
Sean $\left(\Omega_{i},\mathcal{F}_{i}\right)$, $i=1,2,\ldots,$ espacios medibles y $\Omega=\prod_{i=1}^{\infty}\Omega_{i}$ el conjunto de todas las sucesiones $\left(\omega_{1},\omega_{2},\ldots,\right)$ tales que $\omega_{i}\in\Omega_{i}$, $i=1,2,\ldots,$. Si $B^{n}\subset\prod_{i=1}^{\infty}\Omega_{i}$, definimos $B_{n}=\left\{\omega\in\Omega:\left(\omega_{1},\omega_{2},\ldots,\omega_{n}\right)\in B^{n}\right\}$. Al conjunto $B_{n}$ se le llama {\em cilindro} con base $B^{n}$, el cilindro es llamado medible si $B^{n}\in\prod_{i=1}^{\infty}\mathcal{F}_{i}$.
\end{Def}

\begin{Def}\label{Def.Proc.Adaptado}[TSP, Ash \cite{RBA}]
Sea $X\left(t\right),t\geq0$ proceso estoc\'astico, el proceso es adaptado a la familia de $\sigma$-\'algebras $\mathcal{F}_{t}$, para $t\geq0$, si para $s<t$ implica que $\mathcal{F}_{s}\subset\mathcal{F}_{t}$, y $X\left(t\right)$ es $\mathcal{F}_{t}$-medible para cada $t$. Si no se especifica $\mathcal{F}_{t}$ entonces se toma $\mathcal{F}_{t}$ como $\mathcal{F}\left(X\left(s\right),s\leq t\right)$, la m\'as peque\~na $\sigma$-\'algebra de subconjuntos de $\Omega$ que hace que cada $X\left(s\right)$, con $s\leq t$ sea Borel medible.
\end{Def}

\begin{Def}\label{Def.Tiempo.Paro}[TSP, Ash \cite{RBA}]\index{Tiempos de Paro}
Sea $\left\{\mathcal{F}\left(t\right),t\geq0\right\}$ familia creciente de sub $\sigma$-\'algebras. es decir, $\mathcal{F}\left(s\right)\subset\mathcal{F}\left(t\right)$ para $s\leq t$. Un tiempo de paro para $\mathcal{F}\left(t\right)$ es una funci\'on $T:\Omega\rightarrow\left[0,\infty\right]$ tal que $\left\{T\leq t\right\}\in\mathcal{F}\left(t\right)$ para cada $t\geq0$. Un tiempo de paro para el proceso estoc\'astico $X\left(t\right),t\geq0$ es un tiempo de paro para las $\sigma$-\'algebras $\mathcal{F}\left(t\right)=\mathcal{F}\left(X\left(s\right)\right)$.
\end{Def}

\begin{Def}\index{Proceso Adaptado}
Sea $X\left(t\right),t\geq0$ proceso estoc\'astico, con $\left(S,\chi\right)$ espacio de estados. Se dice que el proceso es adaptado a $\left\{\mathcal{F}\left(t\right)\right\}$, es decir, si para cualquier $s,t\in I$, $I$ conjunto de \'indices, $s<t$, se tiene que $\mathcal{F}\left(s\right)\subset\mathcal{F}\left(t\right)$, y $X\left(t\right)$ es $\mathcal{F}\left(t\right)$-medible,
\end{Def}

\begin{Def}\index{Proceso de Markov}
Sea $X\left(t\right),t\geq0$ proceso estoc\'astico, se dice que es un Proceso de Markov relativo a $\mathcal{F}\left(t\right)$ o que $\left\{X\left(t\right),\mathcal{F}\left(t\right)\right\}$ es de Markov si y s\'olo si para cualquier conjunto $B\in\chi$,  y $s,t\in I$, $s<t$ se cumple que
\begin{equation}\label{Prop.Markov}
P\left\{X\left(t\right)\in B|\mathcal{F}\left(s\right)\right\}=P\left\{X\left(t\right)\in B|X\left(s\right)\right\}.
\end{equation}
\end{Def}

\begin{Note}
Si se dice que $\left\{X\left(t\right)\right\}$ es un Proceso de Markov sin mencionar $\mathcal{F}\left(t\right)$, se asumir\'a que 
\begin{eqnarray*}
\mathcal{F}\left(t\right)=\mathcal{F}_{0}\left(t\right)=\mathcal{F}\left(X\left(r\right),r\leq t\right),
\end{eqnarray*}
entonces la ecuaci\'on (\ref{Prop.Markov}) se puede escribir como
\begin{equation}
P\left\{X\left(t\right)\in B|X\left(r\right),r\leq s\right\} = P\left\{X\left(t\right)\in B|X\left(s\right)\right\}.
\end{equation}
\end{Note}

\begin{Teo}
Sea $\left(X_{n},\mathcal{F}_{n},n=0,1,\ldots,\right\}$ Proceso de Markov con espacio de estados $\left(S_{0},\chi_{0}\right)$ generado por una distribuici\'on inicial $P_{o}$ y probabilidad de transici\'on $p_{mn}$, para $m,n=0,1,\ldots,$ $m<n$, que por notaci\'on se escribir\'a como $p\left(m,n,x,B\right)\rightarrow p_{mn}\left(x,B\right)$. Sea $S$ tiempo de paro relativo a la $\sigma$-\'algebra $\mathcal{F}_{n}$. Sea $T$ funci\'on medible, $T:\Omega\rightarrow\left\{0,1,\ldots,\right\}$. Sup\'ongase que $T\geq S$, entonces $T$ es tiempo de paro. Si $B\in\chi_{0}$,
entonces
\begin{equation}\label{Prop.Fuerte.Markov}
P\left\{X\left(T\right)\in B,T<\infty|\mathcal{F}\left(S\right)\right\} = p\left(S,T,X\left(s\right),B\right).
\end{equation}
en $\left\{T<\infty\right\}$.
\end{Teo}

\subsection{Cadenas de Markov}

\begin{Def}\index{Cadena de Markov}
Sea $\left(\Omega,\mathcal{F},\prob\right)$ un espacio de probabilidad y $\mathbf{E}$ un conjunto no vac\'io, finito o numerable. Una sucesi\'on de variables aleatorias $\left\{X_{n}:\Omega\rightarrow\mathbf{E},n\geq0\right\}$ se le llama \textit{Cadena de Markov} con espacio de estados $\mathbf{E}$ si satisface la condici\'on de Markov, esto es, si para todo $n\geq1$ y toda sucesi\'on $x_{0},x_{1},\ldots,x_{n},x,y\in\mathbf{E}$ se cumple que 

\begin{equation}
P\left\{X_{n}=y|X_{n-1}=x,\ldots,X_{0}=x_{0}\right\}=P\left\{X_{n}=x_{n}|X_{n-1}=x_{n-1}\right\}.
\end{equation}
La distribuci\'on de $X_{0}$ se llama distribuci\'on inicial y se denotar\'a por $\pi$.
\end{Def}

\begin{Note}\index{Probabilidades Condicionales}
Las probabilidades condicionales $P\left\{X_{n}=y|X_{n-1}=x\right\}$ se les llama \textit{probabilidades condicionales}\index{Probabilidades Condicionales}
\end{Note}

\begin{Note}\index{Cadenas Homog\'eneas}
En este trabajo se considerar\'an solamente aquellas cadenas de Markov con probabilidades de transici\'on estacionarias, es decir, aquellas que no dependen del valor de $n$ (se dice que es una cadena homog\'enea), es decir, cuando se diga $X_{n},n\geq0$ es cadena de Markov, se entiende que es una sucesi\'on de variables aleatorias que satisfacen la propiedad de Markov y que tienen probabilidades de transici\'on estacionarias.\index{Cadena Homog\'enea}
\end{Note}

\begin{Note}
Para una cadena de Markov Homog\'enea se tiene la siguiente denotaci\'on
\begin{equation}
P\left\{X_{n}=y|X_{n-1}=x\right\}=P_{x,y}.
\end{equation}
\end{Note}

\begin{Note}\index{Probabilidades de Transici\'on}
Para $m\geq1$ se denotar\'a por $P^{(m)}_{x,y}$ a $P\left\{X_{n+m}=y|X_{n}=x\right\}$, que significa la probabilidad de ir en $m$ pasos o unidades de tiempo de $x$ a $y$, y se le llama \textit{probabilidad de transici\'on en $m$ pasos}.
\end{Note}

\begin{Note}\index{Delta de Kronecker}
Para $x,y\in\mathbf{E}$ se define a $P^{(0)}_{x,y}$ como $\delta_{x,y}$, donde $\delta_{x,y}$ es la delta de Kronecker, es decir, vale 1 si $x=y$ y 0 en otro caso.
\end{Note}

\begin{Note}\index{Matriz de Transici\'on}
En el caso de que $\mathbf{E}$ sea finito, se considera la matrix $P=\left(P_{x,y}\right)_{x,y\in \mathbf{E}}$ y se le llama \textit{matriz de transici\'on}.
\end{Note}

\begin{Note}
Si la distribuci\'on inicial $\pi$ es igual al vector $\left(\delta_{x,y}\right)_{y\in\mathbf{E}}$, es decir,
\begin{eqnarray*}
P\left(X_{0}=x)=1\right)\textrm{ y }P\left(X_{0}\neq x\right)=0,
\end{eqnarray*}
entonces se toma la notaci\'on 
\begin{eqnarray}
&&P_{x}\left(A\right)=P\left(A|X_{0}=x\right),A\in\mathcal{F},
\end{eqnarray}
y se dice que la cadena empieza en $A$. Se puede demostrar que $P_{x}$ es una nueva medida de probabilidad en el espacio $\left(\Omega,\mathcal{F}\right)$.
\end{Note}

\begin{Note}
La suma de las entradas de los renglones de la matriz de transici\'on es igual a uno, es decir, para todo $x\in \mathbf{E}$ se tiene $\sum_{y\in\mathbf{E}}P_{x,y}=1$.
\end{Note}

Para poder obtener uno de los resultados m\'as importantes en cadenas de Markov, la \textit{ecuaci\'on de Chapman-kolmogorov} se requieren los siguientes resultados:

\begin{Lema}
Sean $x,y,z\in\Eb$ y $0\leq m\leq n-1$, entonces se cumple que
\begin{equation}
P\left(X_{n+1}=y|X_{n}=z,X_{m}=x\right)=P_{z,y}.
\end{equation}
\end{Lema}

\begin{Prop}
Si $x_{0},x_{1},\ldots,x_{n}\in \Eb$ y $\pi\left(x_{0}\right)=P\left(X_{0}=x_{0}\right)$, entonces
\begin{equation}
P\left(X_{1}=x_{1},\ldots,X_{n}=x_{n},X_{0}=x_{0}\right)=\pi\left(x_{0}\right)P_{x_{0},x_{1}}\cdot P_{x_{1},x_{2}}\cdots P_{x_{n-1},x_{n}}.
\end{equation}
\end{Prop}

De la proposici\'on anterior se tiene
\begin{equation}
P\left(X_{1}=x_{1},\ldots,X_{n}=x_{n}|X_{0}=x_{0}\right)=P_{x_{0},x_{1}}\cdot P_{x_{1},x_{2}}\cdots P_{x_{n-1},x_{n}}.
\end{equation}

finalmente tenemos la siguiente proposici\'on:

\begin{Prop}
Sean $n,k\in\nat$ fijos y $x_{0},x_{1},\ldots,x_{n},\ldots,x_{n+k}\in\Eb$, entonces
\begin{eqnarray*}
&&P\left(X_{n+1}=x_{n+1},\ldots,X_{n+k}=x_{n+k}|X_{n}=x_{n},\ldots,X_{0}=x_{0}\right)\\
&=&P\left(X_{1}=x_{n+1},X_{2}=x_{n+2},\cdots,X_{k}=x_{n+k}|X_{0}=x_{n}\right).
\end{eqnarray*}
\end{Prop}

\begin{Ejem}
Sea $X_{n}$ una variable aleatoria al tiempo $n$ tal que
\begin{eqnarray}
\begin{array}{l}
P\left(X_{n+1}=1 \mid X_{n}=0\right)=p,\\
P\left(X_{n+1}=0 \mid X_{n}=1\right)=q=1-p,\\
P\left(X_{0}=0\right)=\pi_{0}\left(0\right).
\end{array}
\end{eqnarray}

\end{Ejem}

Se puede demostrar que
\begin{eqnarray}
\begin{array}{l}
P\left(X_{n}=0\right)=\frac{q}{p+q},\\
P\left(X_{n}=1\right)=\frac{p}{p+q}.
\end{array}
\end{eqnarray}

\begin{Ejem}
El problema de la Caminata Aleatoria.
\end{Ejem}

\begin{Ejem}
El problema de la ruina del jugador.
\end{Ejem}

\begin{Ejem}
Sea $\left\{Y_{i}\right\}_{i=0}^{\infty}$ sucesi\'on de variables aleatorias independientes e identicamente distribuidas, definidas sobre un espacio de probabilidad $\left(\Omega,\mathcal{F},\prob\right)$ y que toman valores enteros, se tiene que la sucesi\'on $\left\{X_{i}\right\}_{i=0}^{\infty}$ definida por $X_{j}=\sum_{i=0}^{j}Y_{i}$ es una cadena de Markov en el conjunto de los n\'umeros enteros.
\end{Ejem}

\begin{Prop}\index{Ecuaciones de Chapman-Kolmogorov}
Para una cadena de Markov $\left(X_{n}\right)_{n\in\nat}$ con espacio de estados $\Eb$ y para todo $n,m\in \nat$ y toda pareja $x,y\in\Eb$ se cumple
\begin{equation}
P\left(X_{n+m}=y|X_{0}=x\right)=\sum_{z\in\Eb}P_{x,z}^{(m)}P_{z,y}^{(n)}=P_{x,y}^{(n+m)}.
\end{equation}
\end{Prop}

\begin{Note}
Para una cadena de Markov con un n\'umero finito de estados, se puede pensar a $P^{n}$ como la $n$-\'esima potencia de la matriz $P$. Sea $\pi_{0}$ distribuci\'on inicial de la cadena de Markov, como 
\begin{eqnarray}
P\left(X_{n}=y\right)=\sum_{x} P\left(X_{0}=x,X_{n}=y\right)=\sum_{x} P\left(X_{0}=x\right)P\left(X_{n}=y|X_{0}=x\right),
\end{eqnarray}
se puede comprobar que 

\begin{eqnarray}
P\left(X_{n}=y\right)=\sum_{x} \pi_{0}\left(X\right)P^{n}\left(x,y\right).
\end{eqnarray}
\end{Note}

Con lo anterior es posible calcular la distribuici\'on de $X_{n}$ en t\'erminos de la distribuci\'on inicial $\pi_{0}$ y la funci\'on de transici\'on de $n$-pasos $P^{n}$,
\begin{eqnarray}
P\left(X_{n+1}=y\right)=\sum_{x} P\left(X_{n}=x\right)P\left(x,y\right).
\end{eqnarray}
\begin{Note}
Si se conoce la distribuci\'on de $X_{0}$ se puede conocer la distribuci\'on de $X_{1}$.
\end{Note}
\subsection{Procesos de Estados de Markov}

\begin{Teo}
Sea $\left(X_{n},\mathcal{F}_{n},n=0,1,\ldots,\right\}$ Proceso de Markov con espacio de estados $\left(S_{0},\chi_{0}\right)$ generado por una distribuici\'on inicial $P_{o}$ y probabilidad de transici\'on $p_{mn}$, para $m,n=0,1,\ldots,$ $m<n$, que por notaci\'on se escribir\'a como $p\left(m,n,x,B\right)\rightarrow p_{mn}\left(x,B\right)$. Sea $S$ tiempo de paro relativo a la $\sigma$-\'algebra $\mathcal{F}_{n}$. Sea $T$ funci\'on medible, $T:\Omega\rightarrow\left\{0,1,\ldots,\right\}$. Sup\'ongase que $T\geq S$, entonces $T$ es tiempo de paro. Si $B\in\chi_{0}$, entonces
\begin{equation}\label{Prop.Fuerte.Markov2}
P\left\{X\left(T\right)\in B,T<\infty|\mathcal{F}\left(S\right)\right\} =
p\left(S,T,X\left(s\right),B\right),
\end{equation}
en $\left\{T<\infty\right\}$.
\end{Teo}

Sea $K$ conjunto numerable y sea $d:K\rightarrow\nat$ funci\'on. Para $v\in K$, $M_{v}$ es un conjunto abierto de $\rea^{d\left(v\right)}$. Entonces \[E=\bigcup_{v\in K}M_{v}=\left\{\left(v,\zeta\right):v\in K,\zeta\in M_{v}\right\}.\]

Sea $\mathcal{E}$ la clase de conjuntos medibles en $E$:
\[\mathcal{E}=\left\{\bigcup_{v\in K}A_{v}:A_{v}\in \mathcal{M}_{v}\right\}.\] donde $\mathcal{M}$ son los conjuntos de Borel de $M_{v}$. Entonces $\left(E,\mathcal{E}\right)$ es un espacio de Borel. El estado del proceso se denotar\'a por $\mathbf{x}_{t}=\left(v_{t},\zeta_{t}\right)$. La distribuci\'on de $\left(\mathbf{x}_{t}\right)$ est\'a determinada por por los siguientes objetos:

\begin{itemize}
\item[i)] Los campos vectoriales $\left(\mathcal{H}_{v},v\in K\right)$. \item[ii)] Una funci\'on medible $\lambda:E\rightarrow \rea_{+}$. \item[iii)] Una medida de transici\'on $Q:\mathcal{E}\times\left(E\cup\Gamma^{*}\right)\rightarrow\left[0,1\right]$ donde
\begin{equation}
\Gamma^{*}=\bigcup_{v\in K}\partial^{*}M_{v}.
\end{equation}
y
\begin{equation}
\partial^{*}M_{v}=\left\{z\in\partial M_{v}:\mathbf{\mathbf{\phi}_{v}\left(t,\zeta\right)=\mathbf{z}}\textrm{ para alguna }\left(t,\zeta\right)\in\rea_{+}\times M_{v}\right\}.
\end{equation}
donde $\partial M_{v}$ denota  la frontera de $M_{v}$.
\end{itemize}

El campo vectorial $\left(\mathcal{H}_{v},v\in K\right)$ se supone tal que para cada $\mathbf{z}\in M_{v}$ existe una \'unica curva integral $\mathbf{\phi}_{v}\left(t,\zeta\right)$ que satisface la ecuaci\'on

\begin{equation}
\frac{d}{dt}f\left(\zeta_{t}\right)=\mathcal{H}f\left(\zeta_{t}\right),
\end{equation}
con $\zeta_{0}=\mathbf{z}$, para cualquier funci\'on suave $f:\rea^{d}\rightarrow\rea$ y $\mathcal{H}$ denota el operador diferencial de primer orden, con $\mathcal{H}=\mathcal{H}_{v}$ y $\zeta_{t}=\mathbf{\phi}\left(t,\mathbf{z}\right)$. Adem\'as se supone que $\mathcal{H}_{v}$ es conservativo, es decir, las curvas integrales est\'an definidas para todo $t>0$.

Para $\mathbf{x}=\left(v,\zeta\right)\in E$ se denota \[t^{*}\mathbf{x}=inf\left\{t>0:\mathbf{\phi}_{v}\left(t,\zeta\right)\in\partial^{*}M_{v}\right\}.\]

En lo que respecta a la funci\'on $\lambda$, se supondr\'a que para cada $\left(v,\zeta\right)\in E$ existe un $\epsilon>0$ tal que la funci\'on $s\rightarrow\lambda\left(v,\phi_{v}\left(s,\zeta\right)\right)\in E$ es integrable para $s\in\left[0,\epsilon\right)$. La medida de transici\'on $Q\left(A;\mathbf{x}\right)$ es una funci\'on medible de $\mathbf{x}$ para cada $A\in\mathcal{E}$, definida para $\mathbf{x}\in E\cup\Gamma^{*}$ y es una medida de probabilidad en $\left(E,\mathcal{E}\right)$ para cada $\mathbf{x}\in E$.\\

El movimiento del proceso $\left(\mathbf{x}_{t}\right)$ comenzando en $\mathbf{x}=\left(n,\mathbf{z}\right)\in E$, se puede construir de la siguiente manera, def\'inase la funci\'on $F$ por

\begin{equation}
F\left(t\right)=\left\{\begin{array}{ll}
exp\left(-\int_{0}^{t}\lambda\left(n,\phi_{n}\left(s,\mathbf{z}\right)\right)ds\right), & t<t^{*}\left(\mathbf{x}\right),\\
0, & t\geq t^{*}\left(\mathbf{x}\right).
\end{array}\right.
\end{equation}

Sea $T_{1}$ una variable aleatoria tal que $\prob\left[T_{1}>t\right]=F\left(t\right)$, ahora sea la variable aleatoria $\left(N,Z\right)$ con distribuici\'on $Q\left(\cdot;\phi_{n}\left(T_{1},\mathbf{z}\right)\right)$. La trayectoria de $\left(\mathbf{x}_{t}\right)$ para $t\leq T_{1}$ es
\begin{eqnarray*}
\mathbf{x}_{t}=\left(v_{t},\zeta_{t}\right)=\left\{\begin{array}{ll}
\left(n,\phi_{n}\left(t,\mathbf{z}\right)\right), & t<T_{1},\\
\left(N,\mathbf{Z}\right), & t=t_{1}.
\end{array}\right.
\end{eqnarray*}

Comenzando en $\mathbf{x}_{T_{1}}$ se selecciona el siguiente tiempo de intersalto $T_{2}-T_{1}$ lugar del post-salto $\mathbf{x}_{T_{2}}$ de manera similar y as\'i sucesivamente. Este procedimiento nos da una trayectoria determinista por partes $\mathbf{x}_{t}$ con tiempos de salto $T_{1},T_{2},\ldots$. Bajo las condiciones enunciadas para $\lambda,T_{1}>0$  y $T_{1}-T_{2}>0$ para cada $i$, con probabilidad 1. Se supone que se cumple la siguiente condici\'on.

\begin{Sup}[Supuesto 3.1, Davis \cite{Davis}]\label{Sup3.1.Davis}
Sea $N_{t}:=\sum_{t}\indora_{\left(t\geq t\right)}$ el n\'umero de saltos en $\left[0,t\right]$. Entonces
\begin{equation}
\esp\left[N_{t}\right]<\infty\textrm{ para toda }t.
\end{equation}
\end{Sup}

es un proceso de Markov, m\'as a\'un, es un Proceso Fuerte de Markov, es decir, la Propiedad Fuerte de Markov\footnote{Revisar p\'agina 362, y 364 de Davis \cite{Davis}.} se cumple para cualquier tiempo de paro.\index{Propiedad Fuerte de Markov}

\subsection{Clasificaci\'on de Estados}

\begin{Def}\index{Tiempos de Paro}
Para $A$ conjunto en el espacio de estados, se define un tiempo de paro $T_{A}$ de $A$ como
\begin{equation}
T_{A}=min_{n>0}\left(X_{n}\in A\right).
\end{equation}
\end{Def}

\begin{Note}
Si $X_{n}\notin A$ para toda $n>0$, $T_{A}=\infty$, es decir,  $T_{A}$ es el primer tiempo positivo que la cadena de Markov est\'a en $A$.
\end{Note}

Una vez que se tiene la definici\'on anterior se puede demostrar la siguiente igualdad:

\begin{Prop}
$P^{n}\left(x,y\right)=\sum_{m=1}^{n}P_{x}\left(T_{y}=m\right)P^{n.m}\left(y,x\right), n\geq1$.
\end{Prop}
\medskip

\begin{Def}
En una cadena de Markov $\left(X_{n}\right)_{n\in\nat}$ con espacio de estados $\Eb$, matriz de transici\'on $\left(P_{x,y}\right)_{x,y\in\Eb}$ y para $x,y\in\Eb$,  se dice que
\begin{itemize}
\item[a) ]  De $x$ se accede a $y$ si existe $n\geq0$ tal que $P_{x,y}^{(n)}>0$ y se denota por $\left(x\rightarrow y\right)$.

\item[b) ] $x$ y $y$ se comunican entre s\'i, lo que se denota por $\left(x\leftrightarrow y\right)$, si se cumplen $\left(x\rightarrow y\right)$ y $\left(y\rightarrow x\right)$.

\item[c) ] Un estado $x\in\Eb$ es estado recurrente si $$P\left(X_{n}=x\textrm{ para alg\'un }n\in\nat|X_{0}=x \right)\equiv1.$$ \index{Estados recurrentes}

\item[d) ] Un estado $x\in\Eb$ es estado transitorio si $$P\left(X_{n}=x\textrm{ para alg\'un }n\in\nat|X_{0}=x \right)<1.$$ \index{Estados transitorios}

\item[e) ] Un estado $x\in\Eb$ se llama absorbente si $P_{x,x}\equiv1$.\index{Estados absorbentes}
\end{itemize}
\end{Def}

Se tiene el siguiente resultado:

\begin{Prop}
$x\leftrightarrow y$ es una relaci\'on de equivalencia y da lugar a una partici\'on del espacio de estados $\Eb$.
\end{Prop}

\begin{Def}
Para $E$ espacio de estados
\begin{itemize}

\item[a)  ] Se dice que $C\subset \Eb$ es una clase de comunicaci\'on si cualesquiera dos estados de $C$ se comunic\'an entre s\'i.\index{Clases de Comunicaci\'on}

\item[b)  ] Dado $x\in\Eb$, su clase de comunicaci\'on se denota por: $C\left(x\right)=\left\{y\in\Eb:x\leftrightarrow y\right\}$.

\item[c)  ] Se dice que un conjunto de estados  $C\subset \Eb$ es cerrado si ning\'un estado de $\Eb-C$ puede ser accedido desde un estado de $C$.
\end{itemize}
\end{Def}

\begin{Def}\index{Cadena Irreducible}
Sea $\Eb$ espacio de estados, se dice que la cadena es irreducible si cualquiera de las siguientes condiciones, equivalentes entre s\'i,  se cumplen
\begin{enumerate}
\item[a) ] Desde cualquier estado de $\Eb$ se puede acceder a cualquier otro.

\item[b) ] Todos los estados se comunican entre s\'i.

\item[c) ] $C\left(x\right)=\Eb$ para alg\'un $x\in\Eb$.

\item[d) ] $C\left(x\right)=\Eb$ para todo $x\in\Eb$.

\item[e) ] El \'unico conjunto cerrado es el total.
\end{enumerate}
\end{Def}
Por lo tanto tenemos la siguiente proposici\'on:
\begin{Prop}  Sea $\Eb$ espacio de estados y $T$ tiempo de paro, entonces se tiene que
\begin{enumerate}
\item[a) ] Un estado $x\in\Eb$ es recurrente si y s\'olo si $P\left(T_{x}<\infty|x_{0}=x\right)=1$.

\item[b) ] Un estado $x\in\Eb$ es transitorio si y s\'olo si $P\left(T_{x}<\infty|x_{0}=x\right)<1$.

\item[c) ] Un estado $x\in\Eb$ es absorbente si y s\'olo si $P\left(T_{x}=1|x_{0}=x\right)=1$.

\end{enumerate}
\end{Prop}

\subsection{Procesos de Markov}
En esta secci\'on se har\'an las siguientes consideraciones: $E$ es un espacio m\'etrico separable y la m\'etrica $d$ es compatible con la topolog\'ia.

\begin{Def}
Una medida finita, $\lambda$ en la $\sigma$-\'algebra de Borel de un espacio metrizable $E$ se dice cerrada si
\begin{equation}\label{Eq.A2.3}
\lambda\left(E\right)=sup\left\{\lambda\left(K\right):K\textrm{ es
compacto en }E\right\}.
\end{equation}
\end{Def}

\begin{Def}\index{Espacio de Rad\'on}
$E$ es un espacio de Rad\'on si cada medida finita en $\left(E,\mathcal{B}\left(E\right)\right)$ es regular interior o cerrada ({\em tight}).
\end{Def}

El siguiente teorema nos permite tener una mejor caracterizaci\'on de los espacios de Rad\'on:
\begin{Teo}\label{Tma.A2.2}
Sea $E$ espacio separable metrizable. Entonces $E$ es de Rad\'on si y s\'olo s\'i cada medida finita en $\left(E,\mathcal{B}\left(E\right)\right)$ es cerrada.
\end{Teo}

Sea $E$ espacio de estados, tal que $E$ es un espacio de Rad\'on, $\mathcal{B}\left(E\right)$ $\sigma$-\'algebra de Borel en $E$, que se denotar\'a por $\mathcal{E}$.\\

Sea $\left(X,\mathcal{G},\prob\right)$ espacio de probabilidad, $I\subset\rea$ conjunto de \'indices. Sea $\mathcal{F}_{\leq t}$ la $\sigma$-\'algebra natural definida como $\sigma\left\{f\left(X_{r}\right):r\in I, r\leq t,f\in\mathcal{E}\right\}$. Se considerar\'a una $\sigma$-\'algebra m\'as general, $ \left(\mathcal{G}_{t}\right)$ tal que $\left(X_{t}\right)$ sea $\mathcal{E}$-adaptado.

\begin{Def}
Una familia $\left(P_{s,t}\right)$ de kernels de Markov en $\left(E,\mathcal{E}\right)$ indexada por pares $s,t\in I$, con $s\leq t$ es una funci\'on de transici\'on en $\ER$, si  para todo $r\leq s< t$ en $I$ y todo $x\in E$, $B\in\mathcal{E}$,
\begin{equation}\label{Eq.Kernels}
P_{r,t}\left(x,B\right)=\int_{E}P_{r,s}\left(x,dy\right)P_{s,t}\left(y,B\right)\footnote{Ecuaci\'on de Chapman-Kolmogorov}.
\end{equation}
\end{Def}

Se dice que la funci\'on de transici\'on $\KM$ en $\ER$, es la funci\'on de transici\'on para un proceso $\PE$  con valores en $E$ y que satisface la propiedad de Markov\footnote{\begin{equation}\label{Eq.1.4.S}
\prob\left\{H|\mathcal{G}_{t}\right\}=\prob\left\{H|X_{t}\right\}\textrm{ }H\in p\mathcal{F}_{\geq t}.
\end{equation}} (\ref{Eq.1.4.S}) relativa a $\left(\mathcal{G}_{t}\right)$, si

\begin{equation}\label{Eq.1.6.S}
\prob\left\{f\left(X_{t}\right)|\mathcal{G}_{s}\right\}=P_{s,t}f\left(X_{t}\right)\textrm{ }s\leq t\in I,\textrm{ }f\in b\mathcal{E}.
\end{equation}

\begin{Def}\index{Semigrupo de Transici\'on de Markov}
Una familia $\left(P_{t}\right)_{t\geq0}$ de kernels de Markov en $\ER$ es llamada {\em Semigrupo de Transici\'on de Markov} o {\em Semigrupo de Transici\'on} si
\[P_{t+s}f\left(x\right)=P_{t}\left(P_{s}f\right)\left(x\right),\textrm{ }t,s\geq0,\textrm{ }x\in E\textrm{ }f\in b\mathcal{E}.\]
\end{Def}
\begin{Note}
Si la funci\'on de transici\'on $\KM$ es llamada homog\'enea si $P_{s,t}=P_{t-s}$.
\end{Note}

Un proceso de Markov que satisface la ecuaci\'on (\ref{Eq.1.6.S}) con funci\'on de transici\'on homog\'enea $\left(P_{t}\right)$ tiene la propiedad caracter\'istica:
\begin{equation}\label{Eq.1.8.S}\index{Propiedad Simple de Markov}
\prob\left\{f\left(X_{t+s}\right)|\mathcal{G}_{t}\right\}=P_{s}f\left(X_{t}\right)\textrm{ }t,s\geq0,\textrm{ }f\in b\mathcal{E}.
\end{equation}
La ecuaci\'on anterior es la {\em Propiedad Simple de Markov} de $X$ relativa a $\left(P_{t}\right)$. En este sentido el proceso $\PE$ cumple con la propiedad de Markov (\ref{Eq.1.8.S}) relativa a $\left(\Omega,\mathcal{G},\mathcal{G}_{t},\prob\right)$ con semigrupo de transici\'on $\left(P_{t}\right)$.\index{Propiedad Simple de Markov}

\section{Procesos de Renovaci\'on y Regenerativos}\label{Procesos.Regenerativos.Estacionarios}

\begin{Def}\label{Def.Tn}\index{Procesos de Conteo}
Sean $0\leq T_{1}\leq T_{2}\leq \ldots$ son tiempos aleatorios infinitos en los cuales ocurren ciertos eventos. El n\'umero de tiempos $T_{n}$ en el intervalo $\left[0,t\right)$ es

\begin{eqnarray}
N\left(t\right)=\sum_{n=1}^{\infty}\indora\left(T_{n}\leq t\right)\textrm{, para }t\geq0.
\end{eqnarray}
\end{Def}

Si se consideran los puntos $T_{n}$ como elementos de $\rea_{+}$, y $N\left(t\right)$ es el n\'umero de puntos en $\rea$. El proceso denotado por $\left\{N\left(t\right):t\geq0\right\}$ es un proceso puntual en $\rea_{+}$. Los $T_{n}$ son los tiempos de ocurrencia, el proceso puntual $N\left(t\right)$ es simple si su n\'umero de ocurrencias son distintas: $0<T_{1}<T_{2}<\ldots$ casi seguramente.

\begin{Def}\index{Procesos de Renovaci\'on}
Un proceso puntual $N\left(t\right)$ es un proceso de renovaci\'on si los tiempos de interocurrencia $\xi_{n}=T_{n}-T_{n-1}$, para $n\geq1$, son independientes e id\'enticamente distribuidos con distribuci\'on $F$, donde $F\left(0\right)=0$ y $T_{0}=0$. Los $T_{n}$ son llamados tiempos de renovaci\'on, referente a la independencia o renovaci\'on de la informaci\'on estoc\'astica en estos tiempos. Los $\xi_{n}$ son los tiempos de inter-renovaci\'on, y $N\left(t\right)$ es el n\'umero de renovaciones en el intervalo $\left[0,t\right)$.
\end{Def}

\begin{Note}
Para definir un proceso de renovaci\'on para cualquier contexto, solamente hay que especificar una distribuci\'on $F$, con $F\left(0\right)=0$, para los tiempos de inter-renovaci\'on. La funci\'on $F$ en turno define las otras variables aleatorias. De manera formal, existe un espacio de probabilidad y una sucesi\'on de variables aleatorias $\xi_{1},\xi_{2},\ldots$ definidas en este con distribuci\'on $F$. Entonces las otras cantidades son 
\begin{eqnarray}
T_{n}=\sum_{k=1}^{n}\xi_{k},\textrm{ y }N\left(t\right)=\sum_{n=1}^{\infty}\indora\left(T_{n}\leq t\right)\textrm{, donde }T_{n}\rightarrow\infty,
\end{eqnarray}
casi seguramente por la Ley Fuerte de los Grandes N\'umeros.
\end{Note}

\subsection*{Procesos Regenerativos Estacionarios}

\begin{Def}\index{Procesos Regenerativos}
Un proceso estoc\'astico a tiempo continuo $\left\{V\left(t\right),t\geq0\right\}$ es un proceso regenerativo si existe una sucesi\'on de variables aleatorias independientes e id\'enticamente distribuidas $\left\{X_{1},X_{2},\ldots\right\}$, sucesi\'on de renovaci\'on, tal que para cualquier conjunto de Borel $A$, 
\end{Def}

\begin{eqnarray}
\prob\left\{V\left(t\right)\in A|X_{1}+X_{2}+\cdots+X_{R\left(t\right)}=s,\left\{V\left(\tau\right),\tau<s\right\}\right\}=\prob\left\{V\left(t-s\right)\in A|X_{1}>t-s\right\},
\end{eqnarray}
para todo $0\leq s\leq t$, donde $R\left(t\right)=\max\left\{X_{1}+X_{2}+\cdots+X_{j}\leq t\right\}=$n\'umero de renovaciones que ocurren en $\left[0,t\right]$.

\begin{Def}\index{Procesos Estacionarios}
Se define el proceso estacionario, $\left\{V^{*}\left(t\right),t\geq0\right\}$, para $\left\{V\left(t\right),t\geq0\right\}$ por

\begin{eqnarray}
\prob\left\{V\left(t\right)\in A\right\}=\frac{1}{\esp\left[X\right]}\int_{0}^{\infty}\prob\left\{V\left(t+x\right)\in A|X>x\right\}\left(1-F\left(x\right)\right)dx,
\end{eqnarray} 
para todo $t\geq0$ y todo conjunto de Borel $A$.
\end{Def}

\begin{Def}\index{Modificaci\'on Medible}
Una modificaci\'on medible de un proceso $\left\{V\left(t\right),t\geq0\right\}$, es una versi\'on de este, $\left\{V\left(t,w\right)\right\}$ conjuntamente medible para $t\geq0$ y para $w\in S$, $S$ espacio de estados para $\left\{V\left(t\right),t\geq0\right\}$.
\end{Def}

\begin{Teo}\label{Proceso.Regenerativo.No.Negativo}
Sea $\left\{V\left(t\right),t\geq\right\}$ un proceso regenerativo no negativo con modificaci\'on medible. Sea $\esp\left[X\right]<\infty$. Entonces el proceso estacionario dado por la ecuaci\'on anterior est\'a bien definido y tiene funci\'on de distribuci\'on independiente de $t$, adem\'as
\begin{itemize}
\item[i)] 
\begin{eqnarray*}
\esp\left[V^{*}\left(0\right)\right]&=&\frac{\esp\left[\int_{0}^{X}V\left(s\right)ds\right]}{\esp\left[X\right]}.\end{eqnarray*}
\item[ii)] Si $\esp\left[V^{*}\left(0\right)\right]<\infty$, equivalentemente, si $\esp\left[\int_{0}^{X}V\left(s\right)ds\right]<\infty$,entonces
\begin{eqnarray*}
\frac{\int_{0}^{t}V\left(s\right)ds}{t}\rightarrow\frac{\esp\left[\int_{0}^{X}V\left(s\right)ds\right]}{\esp\left[X\right]}.
\end{eqnarray*}
con probabilidad 1 y en media, cuando $t\rightarrow\infty$.
\end{itemize}
\end{Teo}

\begin{Cor}
Sea $\left\{V\left(t\right),t\geq0\right\}$ un proceso regenerativo no negativo, con modificaci\'on medible. Si $\esp V<\infty$, $V$ es no-aritm\'etica, y para todo $x\geq0$, $P\left\{V\left(t\right)\leq x,C>x\right\}$ es de variaci\'on acotada como funci\'on de $t$ en cada intervalo finito $\left[0,\tau\right]$, entonces $V\left(t\right)$ converge en distribuci\'on  cuando $t\rightarrow\infty$ y $$\esp V=\frac{\esp \int_{0}^{X}V\left(s\right)ds}{\esp X}.$$
Donde $V$ tiene la distribuci\'on l\'imite de $V\left(t\right)$ cuando $t\rightarrow\infty$.

\end{Cor}

Para el caso discreto se tienen resultados similares.

%
\subsection{Teorema Principal de Renovaci\'on}

\begin{Note} Una funci\'on $h:\rea_{+}\rightarrow\rea$ es Directamente Riemann Integrable (DRI) en los siguientes casos:
\begin{itemize}
\item[a)] $h\left(t\right)\geq0$ es decreciente y Riemann Integrable.
\item[b)] $h$ es continua excepto posiblemente en un conjunto de Lebesgue de medida 0, y $|h\left(t\right)|\leq b\left(t\right)$, donde $b$ es DRI.
\end{itemize}
\end{Note}

\begin{Teo}[Teorema Principal de Renovaci\'on]\index{Teorema Principal de Renovaci\'on}
Si $F$ es no aritm\'etica y $h\left(t\right)$ es DRI, entonces

\begin{eqnarray}
lim_{t\rightarrow\infty}U\star h=\frac{1}{\mu}\int_{\rea_{+}}h\left(s\right)ds.
\end{eqnarray}
\end{Teo}

\begin{Prop}
Cualquier funci\'on $H\left(t\right)$ acotada en intervalos finitos y que es 0 para $t<0$, puede expresarse como
\begin{eqnarray}
H\left(t\right)=U\star h\left(t\right)\textrm{,  donde }h\left(t\right)=H\left(t\right)-F\star H\left(t\right).
\end{eqnarray}
\end{Prop}

\begin{Def}\index{Procesos Crudamente Regenerativos}
Un proceso estoc\'astico $X\left(t\right)$ es crudamente regenerativo en un tiempo aleatorio positivo $T$ si
\begin{eqnarray}
\esp\left[X\left(T+t\right)|T\right]=\esp\left[X\left(t\right)\right]\textrm{, para }t\geq0,\end{eqnarray}
y con las esperanzas anteriores finitas.
\end{Def}

\begin{Prop}
Sup\'ongase que $X\left(t\right)$ es un proceso crudamente regenerativo en $T$, que tiene distribuci\'on $F$. Si $\esp\left[X\left(t\right)\right]$ es acotado en intervalos finitos, entonces
\begin{eqnarray}
\esp\left[X\left(t\right)\right]=U\star h\left(t\right)\textrm{,  donde }h\left(t\right)=\esp\left[X\left(t\right)\indora\left(T>t\right)\right].
\end{eqnarray}
\end{Prop}

\begin{Teo}[Regeneraci\'on Cruda]
Sup\'ongase que $X\left(t\right)$ es un proceso con valores positivo crudamente regenerativo en $T$, y def\'inase $M=\sup\left\{|X\left(t\right)|:t\leq T\right\}$. Si $T$ es no aritm\'etico; y $M$ y $MT$ tienen media finita, entonces
\begin{eqnarray}
lim_{t\rightarrow\infty}\esp\left[X\left(t\right)\right]=\frac{1}{\mu}\int_{\rea_{+}}h\left(s\right)ds,
\end{eqnarray}
donde $h\left(t\right)=\esp\left[X\left(t\right)\indora\left(T>t\right)\right]$.
\end{Teo}

\begin{Def}\index{Trayectorias Muestrales}
Para el proceso $\left\{\left(N\left(t\right),X\left(t\right)\right):t\geq0\right\}$, sus trayectoria muestrales en el intervalo de tiempo $\left[T_{n-1},T_{n}\right)$ est\'an descritas por
\begin{eqnarray*}
\zeta_{n}=\left(\xi_{n},\left\{X\left(T_{n-1}+t\right):0\leq t<\xi_{n}\right\}\right).
\end{eqnarray*}
Este $\zeta_{n}$ es el $n$-\'esimo segmento del proceso. El proceso es regenerativo sobre los tiempos $T_{n}$ si sus segmentos $\zeta_{n}$ son independientes e id\'enticamennte distribuidos.
\end{Def}

\begin{Note}
Si $\tilde{X}\left(t\right)$ con espacio de estados $\tilde{S}$ es regenerativo sobre $T_{n}$, entonces $X\left(t\right)=f\left(\tilde{X}\left(t\right)\right)$ tambi\'en es regenerativo sobre $T_{n}$, para cualquier funci\'on $f:\tilde{S}\rightarrow S$.
\end{Note}

\begin{Note}
Los procesos regenerativos son crudamente regenerativos, pero no al rev\'es.
\end{Note}

\begin{Note}\label{Procesos.Regenerativo.Clasico}
Un proceso estoc\'astico a tiempo continuo o discreto es regenerativo si existe un proceso de renovaci\'on  tal que los segmentos del proceso entre tiempos de renovaci\'on sucesivos son i.i.d., es decir, para $\left\{X\left(t\right):t\geq0\right\}$ proceso estoc\'astico a tiempo continuo con espacio de estados $S$, espacio m\'etrico.
\end{Note}

Consid\'erese $\left\{X\left(t\right):t\geq0\right\}$ Proceso Estoc\'astico a tiempo continuo con estado de espacios $S$, espacio m\'etrico, con trayectorias continuas por la derecha y con l\'imites por la izquierda c.s. Sea $N\left(t\right)$ un proceso de renovaci\'on en $\rea_{+}$ definido en el mismo espacio de probabilidad que $X\left(t\right)$, con tiempos de renovaci\'on $T$ y tiempos de inter-renovaci\'on $\xi_{n}=T_{n}-T_{n-1}$, con misma distribuci\'on $F$ de media finita $\mu$.

\subsection{Propiedades de los Procesos de Renovaci\'on}
Los tiempos $T_{n}$ est\'an relacionados con los conteos de $N\left(t\right)$ por
\begin{eqnarray}
\begin{array}{l}
\left\{N\left(t\right)\geq n\right\}=\left\{T_{n}\leq t\right\},\\
T_{N\left(t\right)}\leq t<T_{N\left(t\right)+1},
\end{array}
\end{eqnarray}

adem\'as $N\left(T_{n}\right)=n$, y 

\begin{eqnarray}
N\left(t\right)=\max\left\{n:T_{n}\leq t\right\}=\min\left\{n:T_{n+1}>t\right\}.
\end{eqnarray}

Por propiedades de la convoluci\'on se sabe que

\begin{eqnarray}
P\left\{T_{n}\leq t\right\}=F^{n\star}\left(t\right),
\end{eqnarray}
que es la $n$-\'esima convoluci\'on de $F$. Entonces 

\begin{eqnarray}
\begin{array}{l}
\left\{N\left(t\right)\geq n\right\}=\left\{T_{n}\leq t\right\},\\
P\left\{N\left(t\right)\leq n\right\}=1-F^{\left(n+1\right)\star}\left(t\right).
\end{array}
\end{eqnarray}

Adem\'as usando el hecho de que $\esp\left[N\left(t\right)\right]=\sum_{n=1}^{\infty}P\left\{N\left(t\right)\geq n\right\}$, se tiene que

\begin{eqnarray}
\esp\left[N\left(t\right)\right]=\sum_{n=1}^{\infty}F^{n\star}\left(t\right).
\end{eqnarray}

\begin{Prop}
Para cada $t\geq0$, la funci\'on generadora de momentos $\esp\left[e^{\alpha N\left(t\right)}\right]$, existe para alguna $\alpha$ en una vecindad del 0, y de aqu\'i que $\esp\left[N\left(t\right)^{m}\right]<\infty$, para $m\geq1$.
\end{Prop}

\begin{Ejem}[\textbf{Proceso Poisson}]\index{Proceso Poisson}
Suponga que se tienen tiempos de inter-renovaci\'on \textit{i.i.d.} del proceso de renovaci\'on $N\left(t\right)$ tienen distribuci\'on exponencial $F\left(t\right)=q-e^{-\lambda t}$ con tasa $\lambda$. Entonces $N\left(t\right)$ es un proceso Poisson con tasa $\lambda$.

\end{Ejem}


\begin{Note}\index{Proceso de Renovaci\'on Retardado}
Si el primer tiempo de renovaci\'on $\xi_{1}$ no tiene la misma distribuci\'on que el resto de las $\xi_{n}$, para $n\geq2$, a $N\left(t\right)$ se le llama Proceso de Renovaci\'on retardado, donde si $\xi$ tiene distribuci\'on $G$, entonces el tiempo $T_{n}$ de la $n$-\'esima renovaci\'on tiene distribuci\'on $G\star F^{\left(n-1\right)\star}\left(t\right)$.
\end{Note}

\begin{Teo}
Para una constante $\mu\leq\infty$ ( o variable aleatoria), las siguientes expresiones son equivalentes:

\begin{eqnarray}
\begin{array}{l}
lim_{n\rightarrow\infty}n^{-1}T_{n}=\mu,\textrm{ c.s.}\\
lim_{t\rightarrow\infty}t^{-1}N\left(t\right)=1/\mu,\textrm{ c.s.}
\end{array}
\end{eqnarray}
\end{Teo}

Es decir, $T_{n}$ satisface la Ley Fuerte de los Grandes N\'umeros s\'i y s\'olo s\'i $N\left(t\right)$ la cumple.

\begin{Cor}[Ley Fuerte de los Grandes N\'umeros para Procesos de Renovaci\'on]\index{Ley Fuerte de los Grandes N\'umeros para Procesos de Renovaci\'on}
Si $N\left(t\right)$ es un proceso de renovaci\'on cuyos tiempos de inter-renovaci\'on tienen media $\mu\leq\infty$, entonces
\begin{eqnarray}
t^{-1}N\left(t\right)\rightarrow 1/\mu,\textrm{ c.s. cuando }t\rightarrow\infty.
\end{eqnarray}
\end{Cor}

Consideremos el proceso estoc\'astico de valores reales $\left\{Z\left(t\right):t\geq0\right\}$ en el mismo espacio de probabilidad que $N\left(t\right)$.

\begin{Def}
Para el proceso $\left\{Z\left(t\right):t\geq0\right\}$, se define la fluctuaci\'on m\'axima de $Z\left(t\right)$ en el intervalo $\left(T_{n-1},T_{n}\right]$:
\begin{eqnarray*}
M_{n}=\sup_{T_{n-1}<t\leq T_{n}}|Z\left(t\right)-Z\left(T_{n-1}\right)|.
\end{eqnarray*}
\end{Def}

\begin{Teo}
Sup\'ongase que $n^{-1}T_{n}\rightarrow\mu$ c.s. cuando $n\rightarrow\infty$, donde $\mu\leq\infty$ es una constante o variable aleatoria. Sea $a$ una constante o variable aleatoria que puede ser infinita cuando $\mu$ es finita, y considere las expresiones l\'imite:
\begin{eqnarray}
\begin{array}{l}
lim_{n\rightarrow\infty}n^{-1}Z\left(T_{n}\right)=a,\textrm{ c.s.}\\
lim_{t\rightarrow\infty}t^{-1}Z\left(t\right)=a/\mu,\textrm{ c.s.}
\end{array}
\end{eqnarray}
La segunda expresi\'on implica la primera. Conversamente, la primera implica la segunda si el proceso $Z\left(t\right)$ es creciente, o si $lim_{n\rightarrow\infty}n^{-1}M_{n}=0$ c.s.
\end{Teo}

\begin{Cor}
Si $N\left(t\right)$ es un proceso de renovaci\'on, y $\left(Z\left(T_{n}\right)-Z\left(T_{n-1}\right),M_{n}\right)$, para $n\geq1$, son variables aleatorias independientes e id\'enticamente distribuidas con media finita, entonces,
\begin{eqnarray}
lim_{t\rightarrow\infty}t^{-1}Z\left(t\right)\rightarrow\frac{\esp\left[Z\left(T_{1}\right)-Z\left(T_{0}\right)\right]}{\esp\left[T_{1}\right]},\textrm{ c.s. cuando  }t\rightarrow\infty.
\end{eqnarray}
\end{Cor}

\subsection{Funci\'on de Renovaci\'on}

Sup\'ongase que $N\left(t\right)$ es un proceso de renovaci\'on con distribuci\'on $F$ con media finita $\mu$.

\begin{Def}\index{Funci\'on de Renovaci\'on}
La funci\'on de renovaci\'on asociada con la distribuci\'on $F$, del proceso $N\left(t\right)$, es
\begin{eqnarray}
U\left(t\right)=\sum_{n=1}^{\infty}F^{n\star}\left(t\right),\textrm{   }t\geq0,
\end{eqnarray}
donde $F^{0\star}\left(t\right)=\indora\left(t\geq0\right)$.
\end{Def}

\begin{Prop}
Sup\'ongase que la distribuci\'on de inter-renovaci\'on $F$ tiene densidad $f$. Entonces $U\left(t\right)$ tambi\'en tiene densidad, para $t>0$, y es $U^{'}\left(t\right)=\sum_{n=0}^{\infty}f^{n\star}\left(t\right)$. Adem\'as
\begin{eqnarray}
\prob\left\{N\left(t\right)>N\left(t-1\right)\right\}=0\textrm{,   }t\geq0.
\end{eqnarray}
\end{Prop}

\begin{Def}\index{Transformada de Laplace-Stieljes}
La Transformada de Laplace-Stieljes de $F$ est\'a dada por
\begin{eqnarray*}
\hat{F}\left(\alpha\right)=\int_{\rea_{+}}e^{-\alpha t}dF\left(t\right)\textrm{,  }\alpha\geq0.
\end{eqnarray*}
\end{Def}

Entonces
\begin{eqnarray}
\hat{U}\left(\alpha\right)=\sum_{n=0}^{\infty}\hat{F^{n\star}}\left(\alpha\right)=\sum_{n=0}^{\infty}\hat{F}\left(\alpha\right)^{n}=\frac{1}{1-\hat{F}\left(\alpha\right)}.
\end{eqnarray}

\begin{Prop}
La Transformada de Laplace $\hat{U}\left(\alpha\right)$ y $\hat{F}\left(\alpha\right)$ determina una a la otra de manera \'unica por la relaci\'on $\hat{U}\left(\alpha\right)=\frac{1}{1-\hat{F}\left(\alpha\right)}$.
\end{Prop}

\begin{Note}\index{Proceso Poisson}
Un proceso de renovaci\'on $N\left(t\right)$ cuyos tiempos de inter-renovaci\'on tienen media finita, es un proceso Poisson con tasa $\lambda$ s\'i y s\'olo s\'i $\esp\left[U\left(t\right)\right]=\lambda t$, para $t\geq0$.
\end{Note}

\begin{Teo}
Sea $N\left(t\right)$ un proceso puntual simple con puntos de localizaci\'on $T_{n}$ tal que $\eta\left(t\right)=\esp\left[N\left(t\right)\right]$ es finita para cada $t$. Entonces para cualquier funci\'on $f:\rea_{+}\rightarrow\rea$,
\begin{eqnarray}
\esp\left[\sum_{n=1}^{N\left(t\right)}f\left(T_{n}\right)\right]=\int_{\left(0,t\right]}f\left(s\right)d\eta\left(s\right)\textrm{,  }t\geq0,
\end{eqnarray}
suponiendo que la integral exista. Adem\'as si $X_{1},X_{2},\ldots$ son variables aleatorias definidas en el mismo espacio de probabilidad que el proceso $N\left(t\right)$ tal que $\esp\left[X_{n}|T_{n}=s\right]=f\left(s\right)$, independiente de $n$. Entonces
\begin{eqnarray}
\esp\left[\sum_{n=1}^{N\left(t\right)}X_{n}\right]=\int_{\left(0,t\right]}f\left(s\right)d\eta\left(s\right)\textrm{,  }t\geq0,
\end{eqnarray} 
suponiendo que la integral exista. 
\end{Teo}

\begin{Cor}[Identidad de Wald para Renovaciones]\index{Identidad de Wald para Renovaciones}
Para el proceso de renovaci\'on $N\left(t\right)$,
\begin{eqnarray}
\esp\left[T_{N\left(t\right)+1}\right]=\mu\esp\left[N\left(t\right)+1\right]\textrm{,  }t\geq0.
\end{eqnarray}  
\end{Cor}

\begin{Def}\index{Ecuaci\'on de Renovaci\'on}
Sea $h\left(t\right)$ funci\'on de valores reales en $\rea$ acotada en intervalos finitos e igual a cero para $t<0$ La ecuaci\'on de renovaci\'on para $h\left(t\right)$ y la distribuci\'on $F$ es
\begin{eqnarray}\label{Ec.Renovacion}
H\left(t\right)=h\left(t\right)+\int_{\left[0,t\right]}H\left(t-s\right)dF\left(s\right)\textrm{,    }t\geq0,
\end{eqnarray}
donde $H\left(t\right)$ es una funci\'on de valores reales. Esto es $H=h+F\star H$. Decimos que $H\left(t\right)$ es soluci\'on de esta ecuaci\'on si satisface la ecuaci\'on, y es acotada en intervalos finitos e iguales a cero para $t<0$.
\end{Def}

\begin{Prop}
La funci\'on $U\star h\left(t\right)$ es la \'unica soluci\'on de la ecuaci\'on de renovaci\'on (\ref{Ec.Renovacion}).
\end{Prop}

\begin{Teo}[Teorema Renovaci\'on Elemental]\index{Teorema Renovaci\'on Elemental}
\begin{eqnarray*}
t^{-1}U\left(t\right)\rightarrow 1/\mu\textrm{,    cuando }t\rightarrow\infty.
\end{eqnarray*}
\end{Teo}


\section{Teor\'ia de Procesos Regenerativos}\label{Thorisson}
\begin{Def}[Definici\'on Cl\'asica]\index{Procesos Regenerativos}
Un proceso estoc\'astico $X=\left\{X\left(t\right):t\geq0\right\}$ es llamado regenerativo si existe una variable aleatoria $R_{1}>0$ tal que
\begin{itemize}
\item[i)] $\left\{X\left(t+R_{1}\right):t\geq0\right\}$ es independiente de $\left\{\left\{X\left(t\right):t<R_{1}\right\},t\geq0\right\}$.
\item[ii)] $\left\{X\left(t+R_{1}\right):t\geq0\right\}$ es estoc\'asticamente equivalente a $\left\{X\left(t\right):t>0\right\}$.
\end{itemize}

Llamamos a $R_{1}$ tiempo de regeneraci\'on, y decimos que $X$ se regenera en este punto.\index{Tiempos de Regeneraci\'on}
\end{Def}

$\left\{X\left(t+R_{1}\right)\right\}$ es regenerativo con tiempo de regeneraci\'on $R_{2}$, independiente de $R_{1}$ pero con la misma distribuci\'on que $R_{1}$. Procediendo de esta manera se obtiene una secuencia de variables aleatorias independientes e id\'enticamente distribuidas $\left\{R_{n}\right\}$ llamados longitudes de ciclo. Si definimos a $Z_{k}\equiv R_{1}+R_{2}+\cdots+R_{k}$, se tiene un proceso de renovaci\'on llamado proceso de renovaci\'on encajado para $X$.\index{Procesos de Renovaci\'on Encajados}

\begin{Note}
La existencia de un primer tiempo de regeneraci\'on, $R_{1}$, implica la existencia de una sucesi\'on completa de estos tiempos $R_{1},R_{2}\ldots,$ que satisfacen la propiedad deseada \cite{Sigman2}.
\end{Note}

\begin{Note} Para la cola $GI/GI/1$ los usuarios arriban con tiempos $t_{n}$ y son atendidos con tiempos de servicio $S_{n}$, los tiempos de arribo forman un proceso de renovaci\'on  con tiempos entre arribos independientes e identicamente distribuidos (\texttt{i.i.d.})$T_{n}=t_{n}-t_{n-1}$, adem\'as los tiempos de servicio son \texttt{i.i.d.} e independientes de los procesos de arribo. Por \textit{estable} se entiende que $\esp S_{n}<\esp T_{n}<\infty$.
\end{Note}

\begin{Def}
Para $x$ fijo y para cada $t\geq0$, sea $I_{x}\left(t\right)=1$ si $X\left(t\right)\leq x$,  $I_{x}\left(t\right)=0$ en caso contrario, y def\'inanse los tiempos promedio

\begin{eqnarray}
\begin{array}{l}
\overline{X}=lim_{t\rightarrow\infty}\frac{1}{t}\int_{0}^{\infty}X\left(u\right)du,\\
\prob\left(X_{\infty}\leq x\right)=lim_{t\rightarrow\infty}\frac{1}{t}\int_{0}^{\infty}I_{x}\left(u\right)du,
\end{array}
\end{eqnarray}

cuando estos l\'imites existan.
\end{Def}

Como consecuencia del teorema de Renovaci\'on-Recompensa, se tiene que el primer l\'imite  existe y es igual a la constante
\begin{eqnarray}
\overline{X}&=&\frac{\esp\left[\int_{0}^{R_{1}}X\left(t\right)dt\right]}{\esp\left[R_{1}\right]},
\end{eqnarray}
suponiendo que ambas esperanzas son finitas.

\begin{Note}
Funciones de procesos regenerativos son regenerativas, es decir, si $X\left(t\right)$ es regenerativo y se define el proceso $Y\left(t\right)$ por $Y\left(t\right)=f\left(X\left(t\right)\right)$ para alguna funci\'on Borel medible $f\left(\cdot\right)$. Adem\'as $Y$ es regenerativo con los mismos tiempos de renovaci\'on que $X$. \\

En general, los tiempos de renovaci\'on, $Z_{k}$ de un proceso regenerativo no requieren ser tiempos de paro con respecto a la evoluci\'on de $X\left(t\right)$.
\end{Note} 

\begin{Note}
Una funci\'on de un proceso de Markov, usualmente no ser\'a un proceso de Markov, sin embargo ser\'a regenerativo si el proceso de Markov lo es.
\end{Note}

\begin{Note}\index{Procesos Positivo Recurrente}
Un proceso regenerativo con media de la longitud de ciclo finita es llamado positivo recurrente.
\end{Note}

\begin{Note}
\begin{itemize}
\item[a)] Si el proceso regenerativo $X$ es positivo recurrente y tiene trayectorias muestrales no negativas, entonces la ecuaci\'on anterior es v\'alida.
\item[b)] Si $X$ es positivo recurrente regenerativo, podemos construir una \'unica versi\'on estacionaria de este proceso, $X_{e}=\left\{X_{e}\left(t\right)\right\}$, donde $X_{e}$ es un proceso estoc\'astico regenerativo y estrictamente estacionario, con distribuci\'on marginal distribuida como $X_{\infty}$
\end{itemize}
\end{Note}

%
\subsection{Procesos de Renovaci\'on y Regenerativos}
%

\begin{Def}[Renewal Process Trinity]
Para un proceso de renovaci\'on $N\left(t\right)$, los siguientes procesos proveen de informaci\'on sobre los tiempos de renovaci\'on.
\begin{itemize}
\item $A\left(t\right)=t-T_{N\left(t\right)}$, el tiempo de recurrencia hacia atr\'as al tiempo $t$, que es el tiempo desde la \'ultima renovaci\'on para $t$.

\item $B\left(t\right)=T_{N\left(t\right)+1}-t$, el tiempo de recurrencia hacia adelante al tiempo $t$, residual del tiempo de renovaci\'on, que es el tiempo para la pr\'oxima renovaci\'on despu\'es de $t$.

\item $L\left(t\right)=\xi_{N\left(t\right)+1}=A\left(t\right)+B\left(t\right)$, la longitud del intervalo de renovaci\'on que contiene a $t$.
\end{itemize}
\end{Def}

\begin{Note}\index{Procesos de Markov}
El proceso tridimensional $\left(A\left(t\right),B\left(t\right),L\left(t\right)\right)$ es regenerativo sobre $T_{n}$, y por ende cada proceso lo es. Cada proceso $A\left(t\right)$ y $B\left(t\right)$ son procesos de Markov a tiempo continuo con trayectorias continuas por partes en el espacio de estados $\rea_{+}$. Una expresi\'on conveniente para su distribuci\'on conjunta es, para $0\leq x<t,y\geq0$ tal que,
\begin{equation}\label{NoRenovacion}
P\left\{A\left(t\right)>x,B\left(t\right)>y\right\}=
P\left\{N\left(t+y\right)-N\left((t-x)\right)=0\right\}.
\end{equation}
\end{Note}

\begin{Ejem}[Tiempos de recurrencia Poisson]
Si $N\left(t\right)$ es un proceso Poisson con tasa $\lambda$, entonces de la expresi\'on (\ref{NoRenovacion}) se tiene que
\begin{eqnarray*}
\begin{array}{l}
P\left\{A\left(t\right)>x,B\left(t\right)>y\right\}=e^{-\lambda\left(x+y\right)}\textrm{, }0\leq x<t,y\geq0,
\end{array}
\end{eqnarray*}
que es la probabilidad Poisson de no renovaciones en un intervalo de longitud $x+y$.
\end{Ejem}

\begin{Note}\index{Cadena de Markov Erg\'odica}
Una cadena de Markov erg\'odica tiene la propiedad de ser estacionaria si la distribuci\'on de su estado al tiempo $0$ es su distribuci\'on estacionaria.
\end{Note}

\begin{Def}\index{Procesos Estacionarios}
Un proceso estoc\'astico a tiempo continuo $\left\{X\left(t\right):t\geq0\right\}$ en un espacio general es estacionario si sus distribuciones finito dimensionales son invariantes bajo cualquier  traslado: para cada $0\leq s_{1}<s_{2}<\cdots<s_{k}$ y $t\geq0$,
\begin{eqnarray*}
\left(X\left(s_{1}+t\right),\ldots,X\left(s_{k}+t\right)\right)=_{d}\left(X\left(s_{1}\right),\ldots,X\left(s_{k}\right)\right).
\end{eqnarray*}
\end{Def}

\begin{Note}
Un proceso de Markov es estacionario si $X\left(t\right)=_{d}X\left(0\right)$, $t\geq0$.
\end{Note}

Considerese el proceso $N\left(t\right)=\sum_{n}\indora\left(\tau_{n}\leq t\right)$ en $\rea_{+}$, con puntos $0<\tau_{1}<\tau_{2}<\cdots$.

\begin{Prop}
Si $N$ es un proceso puntual estacionario y $\esp\left[N\left(1\right)\right]<\infty$, entonces 
\begin{eqnarray}
\esp\left[N\left(t\right)\right]=t\esp\left[N\left(1\right)\right]\textrm{, }t\geq0.
\end{eqnarray}
\end{Prop}

\begin{Teo}
Los siguientes enunciados son equivalentes
\begin{itemize}
\item[i)] El proceso retardado de renovaci\'on $N$ es estacionario,
\item[ii)] EL proceso de tiempos de recurrencia hacia adelante $B\left(t\right)$ es estacionario,
\item[iii)] $\esp\left[N\left(t\right)\right]=t/\mu$,
\item[iv)] $G\left(t\right)=F_{e}\left(t\right)=\frac{1}{\mu}\int_{0}^{t}\left[1-F\left(s\right)\right]ds$.
\end{itemize}
Cuando estos enunciados son ciertos, $P\left\{B\left(t\right)\leq x\right\}=F_{e}\left(x\right)$, para $t,x\geq0$.
\end{Teo}

\begin{Note}
Una consecuencia del teorema anterior es que el Proceso Poisson es el \'unico proceso sin retardo que es estacionario.
\end{Note}

\begin{Cor}
El proceso de renovaci\'on $N\left(t\right)$ sin retardo, y cuyos tiempos de inter renovaci\'on tienen media finita, es estacionario si y s\'olo si es un proceso Poisson.
\end{Cor}

\subsection{Procesos Regenerativos Estacionarios}

Para $\left\{X\left(t\right):t\geq0\right\}$ Proceso Estoc\'astico a tiempo continuo con estado de espacios $S$, que es un espacio m\'etrico, con trayectorias continuas por la derecha y con l\'imites por la izquierda c.s. Sea $N\left(t\right)$ un proceso de renovaci\'on en $\rea_{+}$ definido en el mismo espacio de probabilidad que $X\left(t\right)$, con tiempos de renovaci\'on $T$ y tiempos de inter-renovaci\'on $\xi_{n}=T_{n}-T_{n-1}$, con misma distribuci\'on $F$ de media finita $\mu$.

\begin{Def}
Un elemento aleatorio en un espacio medible $\left(E,\mathcal{E}\right)$ en un espacio de probabilidad $\left(\Omega,\mathcal{F},\prob\right)$ a $\left(E,\mathcal{E}\right)$, es decir, para $A\in \mathcal{E}$,  se tiene que $\left\{Y\in A\right\}\in\mathcal{F}$, donde
 \begin{eqnarray}
 \left\{Y\in A\right\}:=\left\{w\in\Omega:Y\left(w\right)\in A\right\}=:Y^{-1}A.
 \end{eqnarray}
\end{Def}

\begin{Note}
Tambi\'en se dice que $Y$ est\'a soportado por el espacio de probabilidad $\left(\Omega,\mathcal{F},\prob\right)$ y que $Y$ es un mapeo medible de $\Omega$ en $E$, es decir, es $\mathcal{F}/\mathcal{E}$ medible.
\end{Note}

\begin{Def}\index{Espacio Producto}
Para cada $i\in \mathbb{I}$ sea $P_{i}$ una medida de probabilidad en un espacio medible $\left(E_{i},\mathcal{E}_{i}\right)$. Se define el espacio producto $\otimes_{i\in\mathbb{I}}\left(E_{i},\mathcal{E}_{i}\right):=\left(\prod_{i\in\mathbb{I}}E_{i},\otimes_{i\in\mathbb{I}}\mathcal{E}_{i}\right)$, donde $\prod_{i\in\mathbb{I}}E_{i}$ es el producto cartesiano de los $E_{i}$'s, y $\otimes_{i\in\mathbb{I}}\mathcal{E}_{i}$ es la $\sigma$-\'algebra producto, es decir, es la $\sigma$-\'algebra m\'as peque\~na en $\prod_{i\in\mathbb{I}}E_{i}$ que hace al $i$-\'esimo mapeo proyecci\'on en $E_{i}$ medible para toda $i\in\mathbb{I}$ es la $\sigma$-\'algebra inducida por los mapeos proyecci\'on. 
\begin{eqnarray}
\otimes_{i\in\mathbb{I}}\mathcal{E}_{i}:=\sigma\left\{\left\{y:y_{i}\in A\right\}:i\in\mathbb{I}\textrm{ y }A\in\mathcal{E}_{i}\right\}.
\end{eqnarray}
\end{Def}

\begin{Def}
Un espacio de probabilidad $\left(\tilde{\Omega},\tilde{\mathcal{F}},\tilde{\prob}\right)$ es una extensi\'on de otro espacio de probabilidad $\left(\Omega,\mathcal{F},\prob\right)$ si $\left(\tilde{\Omega},\tilde{\mathcal{F}},\tilde{\prob}\right)$ soporta un elemento aleatorio $\xi\in\left(\Omega,\mathcal{F}\right)$ que tienen a $\prob$ como distribuci\'on.
\end{Def}

\begin{Teo}
Sea $\mathbb{I}$ un conjunto de \'indices arbitrario. Para cada $i\in\mathbb{I}$ sea $P_{i}$ una medida de probabilidad en un espacio medible $\left(E_{i},\mathcal{E}_{i}\right)$. Entonces existe una \'unica medida de probabilidad $\otimes_{i\in\mathbb{I}}P_{i}$ en $\otimes_{i\in\mathbb{I}}\left(E_{i},\mathcal{E}_{i}\right)$ tal que 

\begin{eqnarray}
\otimes_{i\in\mathbb{I}}P_{i}\left(y\in\prod_{i\in\mathbb{I}}E_{i}:y_{i}\in A_{i_{1}},\ldots,y_{n}\in A_{i_{n}}\right)=P_{i_{1}}\left(A_{i_{n}}\right)\cdots P_{i_{n}}\left(A_{i_{n}}\right),
\end{eqnarray}
para todos los enteros $n>0$, toda $i_{1},\ldots,i_{n}\in\mathbb{I}$ y todo $A_{i_{1}}\in\mathcal{E}_{i_{1}},\ldots,A_{i_{n}}\in\mathcal{E}_{i_{n}}$.
\end{Teo}

La medida $\otimes_{i\in\mathbb{I}}P_{i}$ es llamada la medida producto y $\otimes_{i\in\mathbb{I}}\left(E_{i},\mathcal{E}_{i},P_{i}\right):=\left(\prod_{i\in\mathbb{I}},E_{i},\otimes_{i\in\mathbb{I}}\mathcal{E}_{i},\otimes_{i\in\mathbb{I}}P_{i}\right)$, es llamado espacio de probabilidad producto.

\begin{Def}\index{Espacio Polaco}
Un espacio medible $\left(E,\mathcal{E}\right)$ es \textit{Polaco} si existe una m\'etrica en $E$ tal que $E$ es completo, es decir cada sucesi\'on de Cauchy converge a un l\'imite en $E$, y \textit{separable}, $E$ tienen un subconjunto denso numerable, y tal que $\mathcal{E}$ es generado por conjuntos abiertos.
\end{Def}

\begin{Def}
Dos espacios medibles $\left(E,\mathcal{E}\right)$ y $\left(G,\mathcal{G}\right)$ son Borel equivalentes \textit{isomorfos} si existe una biyecci\'on $f:E\rightarrow G$ tal que $f$ es $\mathcal{E}/\mathcal{G}$ medible y su inversa $f^{-1}$ es $\mathcal{G}/\mathcal{E}$ medible. La biyecci\'on es una equivalencia de Borel.
\end{Def}

\begin{Def}\index{Espacion Est\'andar}
Un espacio medible  $\left(E,\mathcal{E}\right)$ es un \textit{espacio est\'andar} si es Borel equivalente a $\left(G,\mathcal{G}\right)$, donde $G$ es un subconjunto de Borel de $\left[0,1\right]$ y $\mathcal{G}$ son los subconjuntos de Borel de $G$.
\end{Def}

\begin{Note}
Cualquier espacio Polaco es un espacio est\'andar.
\end{Note}

\begin{Def}
Un proceso estoc\'astico con conjunto de \'indices $\mathbb{I}$ y espacio de estados $\left(E,\mathcal{E}\right)$ es una familia $Z=\left(\mathbb{Z}_{s}\right)_{s\in\mathbb{I}}$ donde $\mathbb{Z}_{s}$ son elementos aleatorios definidos en un espacio de probabilidad com\'un $\left(\Omega,\mathcal{F},\prob\right)$ y todos toman valores en $\left(E,\mathcal{E}\right)$.
\end{Def}

\begin{Def}\index{Procesos Continuos por un lado}
Un proceso estoc\'astico \textit{one-sided contiuous time} (\textbf{PEOSCT}) es un proceso estoc\'astico con conjunto de \'indices $\mathbb{I}=\left[0,\infty\right)$.
\end{Def}

Sea $\left(E^{\mathbb{I}},\mathcal{E}^{\mathbb{I}}\right)$ denota el espacio producto $\left(E^{\mathbb{I}},\mathcal{E}^{\mathbb{I}}\right):=\otimes_{s\in\mathbb{I}}\left(E,\mathcal{E}\right)$. Vamos a considerar $\mathbb{Z}$ como un mapeo aleatorio, es decir, como un elemento aleatorio en $\left(E^{\mathbb{I}},\mathcal{E}^{\mathbb{I}}\right)$ definido por $Z\left(w\right)=\left(Z_{s}\left(w\right)\right)_{s\in\mathbb{I}}$ y $w\in\Omega$.

\begin{Note}
La distribuci\'on de un proceso estoc\'astico $Z$ es la distribuci\'on de $Z$ como un elemento aleatorio en $\left(E^{\mathbb{I}},\mathcal{E}^{\mathbb{I}}\right)$. La distribuci\'on de $Z$ esta determinada de manera \'unica por las distribuciones finito dimensionales.
\end{Note}

\begin{Note}
En particular cuando $Z$ toma valores reales, es decir, $\left(E,\mathcal{E}\right)=\left(\mathbb{R},\mathcal{B}\right)$ las distribuciones finito dimensionales est\'an determinadas por las funciones de distribuci\'on finito dimensionales

\begin{eqnarray}
\prob\left(Z_{t_{1}}\leq x_{1},\ldots,Z_{t_{n}}\leq x_{n}\right),x_{1},\ldots,x_{n}\in\mathbb{R},t_{1},\ldots,t_{n}\in\mathbb{I},n\geq1.
\end{eqnarray}
\end{Note}

\begin{Note}
Para espacios polacos $\left(E,\mathcal{E}\right)$ el Teorema de Consistencia de Kolmogorov asegura que dada una colecci\'on de distribuciones finito dimensionales consistentes, siempre existe un proceso estoc\'astico que posee tales distribuciones finito dimensionales.
\end{Note}

\begin{Def}\index{Trayectorias}
Las trayectorias de $Z$ son las realizaciones $Z\left(w\right)$, para $w\in\Omega$ del mapeo aleatorio $Z$.
\end{Def}

\begin{Note}
Algunas restricciones se imponen sobre las trayectorias, por ejemplo que sean continuas por la derecha, o continuas por la derecha con l\'imites por la izquierda, o de manera m\'as general, se pedir\'a que caigan en alg\'un subconjunto $H$ de $E^{\mathbb{I}}$. En este caso es natural considerar a $Z$ como un elemento aleatorio que no est\'a en $\left(E^{\mathbb{I}},\mathcal{E}^{\mathbb{I}}\right)$ sino en $\left(H,\mathcal{H}\right)$, donde $\mathcal{H}$ es la $\sigma$-\'algebra generada por los mapeos proyecci\'on que toman a $z\in H$ a $z_{t}\in E$ para $t\in\mathbb{I}$. A $\mathcal{H}$ se le conoce como la traza de $H$ en $E^{\mathbb{I}}$, es decir,
\begin{eqnarray}
\mathcal{H}:=E^{\mathbb{I}}\cap H:=\left\{A\cap H:A\in E^{\mathbb{I}}\right\}.
\end{eqnarray}
\end{Note}

\begin{Note}\index{Espacio de Trayectorias}
$Z$ tiene trayectorias con valores en $H$ y cada $Z_{t}$ es un mapeo medible de $\left(\Omega,\mathcal{F}\right)$ a $\left(H,\mathcal{H}\right)$. Cuando se considera un espacio de trayectorias en particular $H$, al espacio $\left(H,\mathcal{H}\right)$ se le llama el espacio de trayectorias de $Z$.
\end{Note}

\begin{Note}
La distribuci\'on del proceso estoc\'astico $Z$ con espacio de trayectorias $\left(H,\mathcal{H}\right)$ es la distribuci\'on de $Z$ como  un elemento aleatorio en $\left(H,\mathcal{H}\right)$. La distribuci\'on, nuevemente, est\'a determinada de manera \'unica por las distribuciones finito dimensionales.
\end{Note}

\begin{Def}
Sea $Z$ un PEOSCT  con espacio de estados $\left(E,\mathcal{E}\right)$ y sea $T$ un tiempo aleatorio en $\left[0,\infty\right)$. Por $Z_{T}$ se entiende el mapeo con valores en $E$ definido en $\Omega$ en la manera obvia:
\begin{eqnarray}
Z_{T}\left(w\right):=Z_{T\left(w\right)}\left(w\right)\textrm{, }w\in\Omega.
\end{eqnarray}
\end{Def}

\begin{Def}\index{Procesos Conjuntamente Medibles}
Un PEOSCT $Z$ es conjuntamente medible (\textbf{CM}) si el mapeo que toma $\left(w,t\right)\in\Omega\times\left[0,\infty\right)$ a $Z_{t}\left(w\right)\in E$ es $\mathcal{F}\otimes\mathcal{B}\left[0,\infty\right)/\mathcal{E}$ medible.
\end{Def}

\begin{Note}
Un PEOSCT-CM implica que el proceso es medible, dado que $Z_{T}$ es una composici\'on  de dos mapeos continuos: el primero que toma $w$ en $\left(w,T\left(w\right)\right)$ es $\mathcal{F}/\mathcal{F}\otimes\mathcal{B}\left[0,\infty\right)$ medible, mientras que el segundo toma $\left(w,T\left(w\right)\right)$ en $Z_{T\left(w\right)}\left(w\right)$ es $\mathcal{F}\otimes\mathcal{B}\left[0,\infty\right)/\mathcal{E}$ medible.
\end{Note}

\begin{Def}\index{Procesos Can\'onicmente Conjuntamente Medibles}
Un PEOSCT con espacio de estados $\left(H,\mathcal{H}\right)$ es can\'onicamente conjuntamente medible (\textbf{CCM}) si el mapeo $\left(z,t\right)\in H\times\left[0,\infty\right)$ en $Z_{t}\in E$ es $\mathcal{H}\otimes\mathcal{B}\left[0,\infty\right)/\mathcal{E}$ medible.
\end{Def}

\begin{Note}
Un PEOSCTCCM implica que el proceso es CM, dado que un PECCM $Z$ es un mapeo de $\Omega\times\left[0,\infty\right)$ a $E$, es la composici\'on de dos mapeos medibles: el primero, toma $\left(w,t\right)$ en $\left(Z\left(w\right),t\right)$ es $\mathcal{F}\otimes\mathcal{B}\left[0,\infty\right)/\mathcal{H}\otimes\mathcal{B}\left[0,\infty\right)$ medible, y el segundo que toma $\left(Z\left(w\right),t\right)$  en $Z_{t}\left(w\right)$ es $\mathcal{H}\otimes\mathcal{B}\left[0,\infty\right)/\mathcal{E}$ medible. Por tanto CCM es una condici\'on m\'as fuerte que CM.
\end{Note}

\begin{Def}\index{Trayectorias Internamente Shift-Invariantes}
Un conjunto de trayectorias $H$ de un PEOSCT $Z$, es internamente shift-invariante (\textbf{ISI}) si 
\begin{eqnarray}
\left\{\left(z_{t+s}\right)_{s\in\left[0,\infty\right)}:z\in H\right\}=H\textrm{, }t\in\left[0,\infty\right).
\end{eqnarray}
\end{Def}

\begin{Def}\index{Mapeo-shift}
Dado un PEOSCTISI, se define el mapeo-shift $\theta_{t}$, $t\in\left[0,\infty\right)$, de $H$ a $H$ por 
\begin{eqnarray}
\theta_{t}z=\left(z_{t+s}\right)_{s\in\left[0,\infty\right)}\textrm{, }z\in H.
\end{eqnarray}
\end{Def}

\begin{Def}\index{Proceso shift-medible}
Se dice que un proceso $Z$ es shift-medible (\textbf{SM}) si $Z$ tiene un conjunto de trayectorias $H$ que es ISI y adem\'as el mapeo que toma $\left(z,t\right)\in H\times\left[0,\infty\right)$ en $\theta_{t}z\in H$ es $\mathcal{H}\otimes\mathcal{B}\left[0,\infty\right)/\mathcal{H}$ medible.
\end{Def}

\begin{Note}
Un proceso estoc\'astico con conjunto de trayectorias $H$ ISI es shift-medible si y s\'olo si es CCM
\end{Note}

\begin{Note}
\begin{itemize}
\item Dado el espacio polaco $\left(E,\mathcal{E}\right)$ se tiene el  conjunto de trayectorias $D_{E}\left[0,\infty\right)$ que es ISI, entonces cumple con ser CCM.

\item Si $G$ es abierto, podemos cubrirlo por bolas abiertas cuay cerradura este contenida en $G$, y como $G$ es segundo numerable como subespacio de $E$, lo podemos cubrir por una cantidad numerable de bolas abiertas.

\end{itemize}
\end{Note}

\begin{Note}
Los procesos estoc\'asticos $Z$ a tiempo discreto con espacio de estados polaco, tambi\'en tiene un espacio de trayectorias polaco y por tanto tiene distribuciones condicionales regulares.
\end{Note}

\begin{Teo}
El producto numerable de espacios polacos es polaco.
\end{Teo}

\begin{Def}\index{Mapeo Medible}
Sea $\left(\Omega,\mathcal{F},\prob\right)$, espacio de probabilidad que soporta al proceso $Z=\left(Z_{s}\right)_{s\in\left[0,\infty\right)}$ y $S=\left(S_{k}\right)_{0}^{\infty}$ donde $Z$ es un PEOSCTM con espacio de estados $\left(E,\mathcal{E}\right)$  y espacio de trayectorias $\left(H,\mathcal{H}\right)$  y adem\'as $S$ es una sucesi\'on de tiempos aleatorios one-sided que satisfacen la condici\'on $0\leq S_{0}<S_{1}<\cdots\rightarrow\infty$. Considerando $S$ como un mapeo medible de $\left(\Omega,\mathcal{F}\right)$ al espacio sucesi\'on $\left(L,\mathcal{L}\right)$, donde 
\begin{eqnarray}
L=\left\{\left(s_{k}\right)_{0}^{\infty}\in\left[0,\infty\right)^{\left\{0,1,\ldots\right\}}:s_{0}<s_{1}<\cdots\rightarrow\infty\right\},
\end{eqnarray}
donde $\mathcal{L}$ son los subconjuntos de Borel de $L$, es decir, $\mathcal{L}=L\cap\mathcal{B}^{\left\{0,1,\ldots\right\}}$. As\'i el par $\left(Z,S\right)$ es un mapeo medible de  $\left(\Omega,\mathcal{F}\right)$ en $\left(H\times L,\mathcal{H}\otimes\mathcal{L}\right)$. El par $\mathcal{H}\otimes\mathcal{L}^{+}$ denotar\'a la clase de todas las funciones medibles de $\left(H\times L,\mathcal{H}\otimes\mathcal{L}\right)$ en $\left(\left[0,\infty\right),\mathcal{B}\left[0,\infty\right)\right)$.
\end{Def}

\begin{Def}
Sea $\theta_{t}$ el mapeo-shift conjunto de $H\times L$ en $H\times L$ dado por
\begin{eqnarray}
\theta_{t}\left(z,\left(s_{k}\right)_{0}^{\infty}\right)=\theta_{t}\left(z,\left(s_{n_{t-}+k}-t\right)_{0}^{\infty}\right),
\end{eqnarray}
donde 
$n_{t-}=inf\left\{n\geq1:s_{n}\geq t\right\}$.
\end{Def}

\begin{Note}
Con la finalidad de poder realizar los shift's sin complicaciones de medibilidad, se supondr\'a que $Z$ es shit-medible, es decir, el conjunto de trayectorias $H$ es invariante bajo shifts del tiempo y el mapeo que toma $\left(z,t\right)\in H\times\left[0,\infty\right)$ en $z_{t}\in E$ es $\mathcal{H}\otimes\mathcal{B}\left[0,\infty\right)/\mathcal{E}$ medible.
\end{Note}

\begin{Def}\index{Procesos Regenerativo Cl\'asico}
Dado un proceso \textbf{PEOSSM} (Proceso Estoc\'astico One Side Shift Medible) $Z$, se dice regenerativo cl\'asico con tiempos de regeneraci\'on $S$ si 
\begin{eqnarray*}
\theta_{S_{n}}\left(Z,S\right)=\left(Z^{0},S^{0}\right)\textrm{, }n\geq0,
\end{eqnarray*}
y adem\'as $\theta_{S_{n}}\left(Z,S\right)$ es independiente de $\left(\left(Z_{s}\right)s\in\left[0,S_{n}\right),S_{0},\ldots,S_{n}\right)$
Si lo anterior se cumple, al par $\left(Z,S\right)$ se le llama regenerativo cl\'asico.
\end{Def}

\begin{Note}\index{Tiempos de Inter-regeneraci\'on}
Si el par $\left(Z,S\right)$ es regenerativo cl\'asico, entonces las longitudes de los ciclos $X_{1},X_{2},\ldots,$ son i.i.d. e independientes de la longitud del retraso $S_{0}$, es decir, $S$ es un proceso de renovaci\'on. Las longitudes de los ciclos tambi\'en son llamados tiempos de inter-regeneraci\'on y tiempos de ocurrencia.
\end{Note}

\begin{Teo}\label{Tma.Reg.Clasico}
Sup\'ongase que el par $\left(Z,S\right)$ es regenerativo cl\'asico con $\esp\left[X_{1}\right]<\infty$. Entonces $\left(Z^{*},S^{*}\right)$ es una versi\'on estacionaria de $\left(Z,S\right)$. Adem\'as, si $X_{1}$ es lattice con span $d$, entonces $\left(Z^{**},S^{**}\right)$ es una versi\'on periodicamente estacionaria de $\left(Z,S\right)$ con periodo $d$.

\end{Teo}

\begin{Def}\index{Variable aleatoria \textit{spread out}}
Una variable aleatoria $X_{1}$ es \textit{spread out} si existe una $n\geq1$ y una  funci\'on $f\in\mathcal{B}^{+}$ tal que $\int_{\rea}f\left(x\right)dx>0$ con $X_{2},X_{3},\ldots,X_{n}$ copias i.i.d  de $X_{1}$, 
\begin{eqnarray}
\prob\left(X_{1}+\cdots+X_{n}\in B\right)\geq\int_{B}f\left(x\right)dx\textrm{para }B\in\mathcal{B}.
\end{eqnarray}
\end{Def}

\begin{Def}\index{Procesos \textit{wide-sense regenerative}}
Dado un proceso estoc\'astico $Z$ se le llama \textit{wide-sense regenerative} (\textbf{WSR}) con tiempos de regeneraci\'on $S$ si $\theta_{S_{n}}\left(Z,S\right)=\left(Z^{0},S^{0}\right)$ para $n\geq0$ en distribuci\'on y $\theta_{S_{n}}\left(Z,S\right)$ es independiente de $\left(S_{0},S_{1},\ldots,S_{n}\right)$ para $n\geq0$.
Se dice que el par $\left(Z,S\right)$ es WSR si lo anterior se cumple.
\end{Def}

\begin{Note}
\begin{itemize}
\item El proceso de trayectorias $\left(\theta_{s}Z\right)_{s\in\left[0,\infty\right)}$ es WSR con tiempos de regeneraci\'on $S$ pero no es regenerativo cl\'asico.
\item Si $Z$ es cualquier proceso estacionario y $S$ es un proceso de renovaci\'on que es independiente de $Z$, entonces $\left(Z,S\right)$ es WSR pero en general no es regenerativo cl\'asico
\end{itemize}
\end{Note}

\begin{Note}
Para cualquier proceso estoc\'astico $Z$, el proceso de trayectorias $\left(\theta_{s}Z\right)_{s\in\left[0,\infty\right)}$ es siempre un proceso de Markov.
\end{Note}

\begin{Teo}
Supongase que el par $\left(Z,S\right)$ es WSR con $\esp\left[X_{1}\right]<\infty$. Entonces $\left(Z^{*},S^{*}\right)$ en el teorema \ref{Tma.Reg.Clasico} es una versi\'on estacionaria de  $\left(Z,S\right)$.
\end{Teo}

\begin{Teo}
Supongase que $\left(Z,S\right)$ es cycle-stationary con $\esp\left[X_{1}\right]<\infty$. Sea $U$ distribuida uniformemente en $\left[0,1\right)$ e independiente de $\left(Z^{0},S^{0}\right)$ y sea $\prob^{*}$ la medida de probabilidad en $\left(\Omega,\prob\right)$ definida por 

\begin{eqnarray}d\prob^{*}=\frac{X_{1}}{\esp\left[X_{1}\right]}d\prob.\end{eqnarray} 

Sea $\left(Z^{*},S^{*}\right)$ con distribuci\'on $\prob^{*}\left(\theta_{UX_{1}}\left(Z^{0},S^{0}\right)\in\cdot\right)$. Entonces $\left(Z^{}*,S^{*}\right)$ es estacionario,
\begin{eqnarray}
\esp\left[f\left(Z^{*},S^{*}\right)\right]=\esp\left[\int_{0}^{X_{1}}f\left(\theta_{s}\left(Z^{0},S^{0}\right)\right)ds\right]/\esp\left[X_{1}\right]
\end{eqnarray}
$f\in\mathcal{H}\otimes\mathcal{L}^{+}$, and $S_{0}^{*}$ es continuo con funci\'on distribuci\'on $G_{\infty}$ definida por 
\begin{eqnarray}
G_{\infty}\left(x\right):=\frac{\esp\left[X_{1}\right]\wedge x}{\esp\left[X_{1}\right]},
\end{eqnarray}
para $x\geq0$ y densidad $\prob\left[X_{1}>x\right]/\esp\left[X_{1}\right]$, con $x\geq0$.
\end{Teo}

\begin{Teo}\label{Teorema.Importante}
Sea $Z$ un Proceso Estoc\'astico un lado shift-medible \textit{one-sided shift-measurable stochastic process}, (PEOSSM), y $S_{0}$ y $S_{1}$ tiempos aleatorios tales que $0\leq S_{0}<S_{1}$ y
\begin{equation}
\theta_{S_{1}}Z=\theta_{S_{0}}Z\textrm{ en distribuci\'on}.
\end{equation}

Entonces el espacio de probabilidad subyacente $\left(\Omega,\mathcal{F},\prob\right)$ puede extenderse para soportar una sucesi\'on de tiempos aleatorios $S$ tales que

\begin{eqnarray}
\begin{array}{l}
\theta_{S_{n}}\left(Z,S\right)=\left(Z^{0},S^{0}\right),n\geq0,\textrm{ en distribuci\'on},\\
\left(Z,S_{0},S_{1}\right)\textrm{ depende de }\left(X_{2},X_{3},\ldots\right)\textrm{ solamente a traves de }\theta_{S_{1}}Z.
\end{array}
\end{eqnarray}
\end{Teo}

\begin{Cor}\label{Tma.Estacionariedad}\index{Teorema de Estacionariedad}
Bajo las condiciones del Teorema anterior , el par $\left(Z,S\right)$ es regenerativo cl\'asico. Si adem\'as se tiene que $\esp\left[X_{1}\right]<\infty$, entonces existe un par $\left(Z^{*},S^{*}\right)$ que es una vesi\'on estacionaria de $\left(Z,S\right)$.
\end{Cor}

\begin{Def}
Los tiempos aleatorios $S_{n}$ dividen $Z$ en 
\begin{itemize}
\item[a)] un retraso $D=\left(Z_{s}\right)_{s\in\left[0,\infty\right)}$,
\item[b)] una sucesi\'on de ciclos $C_{n}=\left(Z_{S_{n-1}+s}\right)_{ s\in\left[0,X_{n}\right)}$, $n\geq1$,
\item[c)] las longitudes de los ciclos $X_{n}=S_{n}-S_{n-1}$, $n\neq1$.
\end{itemize}
\end{Def}

\begin{Note}\index{Procesos Delay-length}
\begin{itemize}
\item[a)] El retraso $D$ y los ciclos $C_{n}$ son procesos estoc\'asticos que se desvanecen en los tiempos aleatorios $S_{0}$ y $X_{n}$ respectivamente.
\item[b)] Las longitudes de los ciclos $X_{1},X_{2},\ldots$ y el retraso de la longitud (\textit{delay-length}) $S_{0}$ son obtenidos por el mismo mapeo medible de sus respectivos ciclos $C_{1},C_{2},\ldots$ y el retraso $D$. 
\item[c)] El par $\left(Z,S\right)$ es un mapeo medible del retraso y de los ciclos y viceversa.\index{Mapeo medible del retraso}
\end{itemize}
\end{Note}

\begin{Def}
$\left(Z,S\right)$ es \textit{zero-delayed} si $S_{0}\equiv0$. Se define el par \textit{zero-delayed} por
\begin{eqnarray}
\left(Z^{0},S^{0}\right):=\theta_{S_{0}}\left(Z,S\right).
\end{eqnarray}
Entonces $S_{0}^{0}\equiv0$ y $S_{0}^{0}\equiv X_{1}^{0}$, mientras que para $n\geq1$ se tiene que $X_{n}^{0}\equiv X_{n}$ y $C_{n}^{0}\equiv C_{n}$.
\end{Def}

\begin{Def}\index{Ciclo Estacionario}
Se le llama al par $\left(Z,S\right)$ \textbf{ciclo-stacionario} si los ciclos forman una sucesi\'on estacionaria, es decir, con $=^{D}$ denota iguales en distribuci\'on:
\begin{eqnarray}
\left(C_{n+1},C_{n+2},\ldots\right)=^{D}\left(C_{1},C_{2},\ldots\right),n\geq0.
\end{eqnarray}
Ciclo-estacionareidad es equivalente a 
\begin{eqnarray}
\theta_{S_{n}}\left(Z,S\right)=^{D}
\left(Z^{0},S^{0}\right),n\geq0,
\end{eqnarray}
donde $\left(C_{n+1},C_{n+2},\ldots\right)$ y $\theta_{S_{n}}\left(Z,S\right)$ son mapeos medibles de cada uno y que no dependen de $n$.
\end{Def}

\begin{Def}
Un par $\left(Z^{*},S^{*}\right)$ es \textbf{estacionario} si $\theta\left(Z^{*},S^{*}\right)=^{D}
\left(Z^{*},S^{*}\right)$, para $t\geq0$.
\end{Def}

%
\subsection{Procesos de Salida y Procesos Regenerativos}\label{Resultados.Salida}
%
En Sigman, Thorison y Wolff \cite{Sigman2} prueban que para la existencia de un una sucesi\'on infinita no decreciente de tiempos de regeneraci\'on $\tau_{1}\leq\tau_{2}\leq\cdots$ en los cuales el proceso se regenera, basta un tiempo de regeneraci\'on $R_{1}$, donde $R_{j}=\tau_{j}-\tau_{j-1}$. Para tal efecto se requiere la existencia de un espacio de probabilidad $\left(\Omega,\mathcal{F},\prob\right)$, y proceso estoc\'astico $\textit{X}=\left\{X\left(t\right):t\geq0\right\}$ con espacio de estados $\left(S,\mathcal{R}\right)$, con $\mathcal{R}$ $\sigma$-\'algebra.

\begin{Prop}
Si existe una variable aleatoria no negativa $R_{1}$ tal que $\theta_{R1}X=_{D}X$, entonces $\left(\Omega,\mathcal{F},\prob\right)$ puede extenderse para soportar una sucesi\'on estacionaria de variables aleatorias $R=\left\{R_{k}:k\geq1\right\}$, tal que para $k\geq1$,
\begin{eqnarray*}
\theta_{k}\left(X,R\right)=_{D}\left(X,R\right).
\end{eqnarray*}

Adem\'as, para $k\geq1$, $\theta_{k}R$ es condicionalmente independiente de $\left(X,R_{1},\ldots,R_{k}\right)$, dado $\theta_{\tau k}X$.

\end{Prop}

A continuaci\'on se enuncian una lista de resultados para sistemas de espera:

\begin{itemize}
\item Doob en 1953 demostr\'o que el estado estacionario de un proceso de partida en un sistema de espera $M/G/\infty$, es Poisson con la misma tasa que el proceso de arribos.

\item Burke en 1968, fue el primero en demostrar que el estado estacionario de un proceso de salida de una cola $M/M/s$ es un proceso Poisson.

\item Disney en 1973 obtuvo el siguiente resultado:

\begin{Teo}
Para el sistema de espera $M/G/1/L$ con disciplina FIFO, el proceso $\textbf{I}$ es un proceso de renovaci\'on si y s\'olo si el proceso denominado longitud de la cola es estacionario y se cumple cualquiera de los siguientes casos:

\begin{itemize}
\item[a)] Los tiempos de servicio son id\'enticamente cero;
\item[b)] $L=0$, para cualquier proceso de servicio $S$;
\item[c)] $L=1$ y $G=D$;
\item[d)] $L=\infty$ y $G=M$.
\end{itemize}
En estos casos, respectivamente, las distribuciones de interpartida $P\left\{T_{n+1}-T_{n}\leq t\right\}$ son

\begin{itemize}
\item[a)] $1-e^{-\lambda t}$, $t\geq0$;
\item[b)] $1-e^{-\lambda t}*F\left(t\right)$, $t\geq0$;
\item[c)] $1-e^{-\lambda t}*\indora_{d}\left(t\right)$, $t\geq0$;
\item[d)] $1-e^{-\lambda t}*F\left(t\right)$, $t\geq0$.
\end{itemize}
\end{Teo}

\item Finch (1959) mostr\'o que para los sistemas $M/G/1/L$, con $1\leq L\leq \infty$ con distribuciones de servicio dos veces diferenciable, solamente el sistema $M/M/1/\infty$ tiene proceso de salida de renovaci\'on estacionario.

\item King (1971) demostr\'o que un sistema de colas estacionario $M/G/1/1$ tiene sus tiempos de interpartida sucesivas $D_{n}$ y $D_{n+1}$ son independientes, si y s\'olo si, $G=D$, en cuyo caso le proceso de salida es de renovaci\'on.

\item Disney (1973) demostr\'o que el \'unico sistema estacionario $M/G/1/L$, que tiene proceso de salida de renovaci\'on  son los sistemas $M/M/1$ y $M/D/1/1$.

\item El siguiente resultado es de Disney y Koning (1985)
\begin{Teo}\index{Proceso Poisson}
En un sistema de espera $M/G/s$, el estado estacionario del proceso de salida es un proceso Poisson para cualquier distribuci\'on de los tiempos de servicio si el sistema tiene cualquiera de las siguientes cuatro propiedades.

\begin{itemize}
\item[a)] $s=\infty$
\item[b)] La disciplina de servicio es de procesador compartido.
\item[c)] La disciplina de servicio es LCFS y preemptive resume, esto se cumple para $L<\infty$
\item[d)] $G=M$.
\end{itemize}

\end{Teo}

\item El siguiente resultado es de Alamatsaz (1983)

\begin{Teo}
En cualquier sistema de colas $GI/G/1/L$ con $1\leq L<\infty$ y distribuci\'on de interarribos $A$ y distribuci\'on de los tiempos de servicio $B$, tal que $A\left(0\right)=0$, $A\left(t\right)\left(1-B\left(t\right)\right)>0$ para alguna $t>0$ y $B\left(t\right)$ para toda $t>0$, es imposible que el proceso de salida estacionario sea de renovaci\'on.
\end{Teo}

\end{itemize}

Estos resultados aparecen en Daley (1968) \cite{Daley68} para $\left\{T_{n}\right\}$ intervalos de inter-arribo, $\left\{D_{n}\right\}$ intervalos de inter-salida y $\left\{S_{n}\right\}$ tiempos de servicio.

\begin{itemize}
\item Si el proceso $\left\{T_{n}\right\}$ es Poisson, el proceso $\left\{D_{n}\right\}$ es no correlacionado si y s\'olo si es un proceso Poisso, lo cual ocurre si y s\'olo si $\left\{S_{n}\right\}$ son exponenciales negativas.

\item Si $\left\{S_{n}\right\}$ son exponenciales negativas, $\left\{D_{n}\right\}$ es un proceso de renovaci\'on  si y s\'olo si es un proceso Poisson, lo cual ocurre si y s\'olo si $\left\{T_{n}\right\}$ es un proceso Poisson.

\item $\esp\left(D_{n}\right)=\esp\left(T_{n}\right)$.

\item Para un sistema de visitas $GI/M/1$ se tiene el siguiente teorema:

\begin{Teo}\label{Tiempos.Interpartida}
En un sistema estacionario $GI/M/1$ los intervalos de interpartida tienen
\begin{eqnarray*}
\esp\left(e^{-\theta D_{n}}\right)&=&\mu\left(\mu+\theta\right)^{-1}\left[\delta\theta
-\mu\left(1-\delta\right)\alpha\left(\theta\right)\right]
\left[\theta-\mu\left(1-\delta\right)^{-1}\right]\\
\alpha\left(\theta\right)&=&\esp\left[e^{-\theta T_{0}}\right]\\
var\left(D_{n}\right)&=&var\left(T_{0}\right)-\left(\tau^{-1}-\delta^{-1}\right)
2\delta\left(\esp\left(S_{0}\right)\right)^{2}\left(1-\delta\right)^{-1}.
\end{eqnarray*}
\end{Teo}

\begin{Teo}\label{Procesos.Salida}
El proceso de salida de un sistema de colas estacionario $GI/M/1$ es un proceso de renovaci\'on si y s\'olo si el proceso de entrada es un proceso Poisson, en cuyo caso el proceso de salida es un proceso Poisson.
\end{Teo}

\begin{Teo}\label{Intervalos.Interpartida}
Los intervalos de interpartida $\left\{D_{n}\right\}$ de un sistema $M/G/1$ estacionario son no correlacionados si y s\'olo si la distribuci\'on de los tiempos de servicio es exponencial negativa, es decir, el sistema es de tipo  $M/M/1$.

\end{Teo}
\end{itemize}

\section{Conclusiones}

Los resultados presentados en este trabajo se enfocan en el estudio de procesos estocásticos regenerativos y estacionarios. En particular, se proporcionaron condiciones que garantizan la regeneración de estos procesos, es decir, las condiciones bajo las cuales existen los tiempos de regeneración. Inicialmente, la idea fue aplicar estos resultados a sistemas de visitas que están conformados por colas. Es natural hacer referencia a procesos de Markov y cadenas de Markov, para los cuales existen resultados análogos que son, en cierto sentido, más fáciles de entender (Secci\'on \ref{Procesos.Estocasticos}).

En la Sección \ref{Procesos.Regenerativos.Estacionarios}, se hace un recuento de resultados para procesos de renovación y sus propiedades. Es importante mencionar la relevancia de los procesos de tipo Poisson debido a la gran cantidad de resultados en la teoría de colas y sus numerosas aplicaciones. Un enfoque alternativo valioso es el presentado en la Sección \ref{Thorisson}, donde, aunque las definiciones y temas se estudian nuevamente, el enfoque y los elementos utilizados difieren, requiriendo nuevas definiciones que ayudan a construir las condiciones para garantizar los tiempos de renovación y, por tanto, la estacionariedad del proceso.

Un resultado importante es el Teorema \ref{Teorema.Importante} y su Corolario \ref{Tma.Estacionariedad}, así como las dos últimas definiciones presentadas al final de la sección. También se recomienda revisar el Teorema \ref{Proceso.Regenerativo.No.Negativo} y la Nota \ref{Procesos.Regenerativo.Clasico} en la Sección \ref{Procesos.Regenerativos.Estacionarios}. Aunque los resultados presentados aquí abundan en el estudio de los procesos regenerativos estacionarios, sin duda se han omitido definiciones y teoremas que podrían ayudar en la comprensión de este tema. A pesar de ser procesos ampliamente estudiados, este trabajo está limitado por la bibliografía revisada, y sería valioso ampliar la lectura a otros autores que ya han realizado estudios en estos temas para lograr un entendimiento más profundo.

En la Subsección \ref{Resultados.Salida}, se listan algunos resultados para sistemas de espera con tiempos de arribo exponenciales, tiempos de servicio generales y un servidor que opera bajo la política de servicio First In, First Out (\textit{FIFO}). Se presentan casos en los que las capacidades de los sistemas de visita son infinitas o limitadas ($L$). Como se mencionó anteriormente, los tiempos de renovación son cruciales, ya que permiten determinar la estacionariedad de un proceso estocástico. En los sistemas de espera, los tiempos de interés son los tiempos de arribo o, equivalentemente, los tiempos entre dos partidas consecutivas. Un resultado importante es el presentado en el Teorema \ref{Tiempos.Interpartida}, junto con los Teoremas \ref{Procesos.Salida} y \ref{Intervalos.Interpartida}.

Finalmente, se sugiere profundizar en el estudio de estos temas, tanto para cadenas de Markov, sistemas de espera, sistemas de visitas y redes de sistemas de visitas, ya que las aplicaciones de estos resultados pueden ser valiosas en áreas como sistemas de transporte y producción. También es recomendable actualizar la literatura con publicaciones recientes y libros de texto actualizados. Además, el desarrollo de métodos numéricos más eficientes para simular procesos regenerativos podría ampliar significativamente el alcance de sus aplicaciones prácticas.

\end{document}